\documentclass[journal,twoside,web]{ieeecolor}
\usepackage{generic}
\usepackage{cite}
\usepackage{amsmath,amssymb,amsfonts}
\usepackage{algorithmic}
\usepackage{graphicx}
\usepackage{algorithm,algorithmic}
\usepackage{hyperref}
\hypersetup{hidelinks}
\usepackage{textcomp}
\def\BibTeX{{\rm B\kern-.05em{\sc i\kern-.025em b}\kern-.08em
    T\kern-.1667em\lower.7ex\hbox{E}\kern-.125emX}}
{J. Sun and M. Cantoni: Energy-gain control of time-varying systems}

\usepackage{bm}
\usepackage{mathrsfs}
\usepackage{bbm}
\usepackage{epstopdf}

\DeclareMathAlphabet{\mathsfit}{OT1}{lmss}{m}{sl}

\newcommand{\trs}{\prime}
\newcommand{\diag}{\mathrm{diag}}

\let\oldforall\forall
\renewcommand{\forall}{\oldforall \, }
\let\oldexist\exists
\renewcommand{\exists}{\oldexist \: }

\newtheorem{assumption}{Assumption}

\newtheorem{remark}{Remark}
\newtheorem{lemma}{Lemma}

\newtheorem{theorem}{Theorem}
\newtheorem{hypothesis}{Hypothesis}

\newcommand\mc[1]{{\color{black} #1}}
\newcommand\mcr[1]{{\color{black} #1}}
\newcommand\lftd[1]{{#1}}

\title{Energy-Gain Control of Time-Varying Systems: Receding Horizon Approximation
\thanks{Funded in part by the Australian Research Council (Grant Number: DP210103272) and a Melbourne Research Scholarship.}}

\author{Jintao Sun and 
Michael Cantoni%~\IEEEmembership{Member,~IEEE}%
 \thanks{J. Sun was a student with the Department of Electrical and Electronic Engineering, The University of Melbourne, VIC 3010, Australia. E-mail: {\tt jintaos@alumni.unimelb.edu.au}}
\thanks{M. Cantoni is with the Department of Electrical and Electronic Engineering, The University of Melbourne, VIC 3010, Australia. Corresponding author. E-mail: {\tt cantoni@unimelb.edu.au}}
}

\begin{document}

\maketitle

\begin{abstract}
Standard formulations of prescribed worst-case disturbance energy-gain control policies for linear time-varying systems depend on all forward model data. In discrete time, this dependence arises through a backward Riccati recursion. This article is about the infinite-horizon $\ell_2$ gain performance of state feedback policies with only finite receding-horizon preview of the model parameters. The proposed synthesis of controllers subject to such a constraint leverages the strict contraction of lifted Riccati operators under uniform controllability and observability. The main approximation result is a sufficient number of preview steps for the incurred performance loss to remain below any set tolerance, relative to the baseline gain bound of the associated infinite-preview controller. Aspects of the result are explored in a numerical example.
\end{abstract}

\begin{IEEEkeywords}
Finite preview, infinite-horizon performance, \mc{non-stationary systems}, Riccati contraction
\end{IEEEkeywords}

\section{Introduction}

\mc{The synthesis of feedback controllers for linear time-varying systems against infinite-horizon quadratic performance criteria has a long history~\cite{kalman1960contributions,athans1966optimal,kwakernaak1972linear,green1994linear,hassibi1999indefinite}. Generally, the standard formulations are somewhat impractical, as the control policies depend on all future time-varying parameters. In the discrete-time setting of this article, the stabilizing solution of a backward Riccati recursion encodes this dependence~\cite{hager1976convergence,anderson1981detectability,bitmead1991riccati,katayama1994h}. 

In periodic settings, including the time-invariant case, availability of the model data across one period suffices to obtain the infinite-horizon solution of the backward recursion via a finite-dimensional algebraic Riccati equation~\cite{caines1970discrete,willems1971least,hench1994numerical,bittanti2009periodic}. In more general time-varying settings, an infinite-horizon lifting, given all model data, also leads to algebraic characterization of the solution~\cite{halanay1994time,peters1997minimum,dullerud1999new}. However, the practical concern of computation remains a challenge unless a finite structure underlies the parameter variation, e.g., periodicity~\cite{dullerud1999new}, or more general switching within a finite set~\cite{essick2014control}. Furthermore, infinite-horizon lifting would be infeasible in supervisory hierarchies involving online roll-out of finite-horizon plans that determine the model data relevant for low-level control.

This article concerns state feedback controller synthesis subject to finite receding-horizon preview of the time-varying model data without preview of the disturbance. The \mcr{goal} is to ensure that the closed-loop infinite-horizon energy gain does not exceed a specified worst-case bound. The proposed approach involves approximation of the standard infinite-preview formulation to within desired tolerance, via a Riccati contraction property that holds under uniform controllability and observability and an associated finite-horizon lifting. Related work on the disturbance-free optimal linear quadratic regulator can be found in~\cite{sun2023receding},~\cite[Ch.~4]{jintaothesis}. } 

\mc{To enable further elaboration of the contribution, the standard formulation of state feedback $\ell_2$ gain controllers with infinite-horizon preview of the model data is recalled in Section~\ref{subsec:probform}. This provides context for Section~\ref{subsec:maincontrib}, where the proposed finite receding-horizon synthesis of such control policies is outlined. Relevant literature is discussed in Section~\ref{subsec:related}, including well-known work on model predictive control, and more recent work on optimal regret, which by contrast focuses on performance loss relative to a policy with non-causal preview of the disturbance. The main technical developments are mapped out in Section~\ref{subsec:map}.}

\subsection{\mcr{Basic notation} and problem formulation}
\label{subsec:probform}
\mc{$\mathbb{N}$ and $\mathbb{R}$ denote the natural and real numbers, $\mathbb{N}_0\!:=\!\{0\}\cup\mathbb{N}$, and 
$\mathbb{R}_{>\theta}\!:=\!\{\vartheta\in\mathbb{R}:\vartheta>\theta\}$. $\mathbb{R}^n$ denotes the space of vectors with $n\in\mathbb{N}$ real co-ordinates, and $\mathbb{R}^{n\times m}$ the space of
matrices with $n\in\mathbb{N}$ rows and $m\in\mathbb{N}$ columns of real entries. The identity matrix is $I_n\in\mathbb{R}^{n\times n}$, all entries of $0_{m,n}\in\mathbb{R}^{m\times n}$ are zero, $X^\prime \in\mathbb{R}^{m\times n}$ is the transpose of $X\in\mathbb{R}^{n\times m}$, and for non-singular $Y\in\mathbb{R}^{n\times n}$, the inverse is $Y^{-1}\in\mathbb{R}^{n\times n}$. Given the sets of symmetric matrices $\mathbb{S}^{n}:=\{Z\in\mathbb{R}^{n\times n}~|~Z=Z^\prime\}$, $\mathbb{S}_{\smash{+}}^{n}:=\{Z\in\mathbb{S}^n~|~(\forall v\in\mathbb{R}^n)~v^\trs Z v\geq 0\}$, and $\mathbb{S}_{\smash{++}}^n:=\{Z\in\mathbb{S}^n~|~(\exists \epsilon>0)(\forall v\in\mathbb{R}^n)~v^\trs Z v \geq \epsilon v^\trs v\}$, $Y \preceq Z$ means $(Z-Y)\in\mathbb{S}_{\smash{+}}^n$, and $Y \prec Z$ means $(Z-Y)\in\mathbb{S}_{\smash{++}}^n$. For $Z\in\mathbb{S}_{\smash{+}}^n$, the matrix square root is $Z^{\smash{1\!/2}}\in\mathbb{S}_{\smash{+}}^n$. For $Y\in\mathbb{S}^n$, $\lambda_{\min}(Y)$ and $\lambda_{\max}(Y)$ denote minimum and maximum eigenvalue.}

Consider the linear time-varying system
\begin{align}\label{eq:ltv_sys_w}
x_{t+1} = A_tx_t + B_tu_t + w_t, \quad t \!\in\! \mathbb{N}_0,
\end{align}
with initial state $x_0=0 \in \mathbb{R}^{n}$, control and disturbance inputs $u_t\in\mathbb{R}^m$ and $w_t\in\mathbb{R}^n$, \mc{respectively, and performance output}
\begin{align}\label{eq:output_z_def}
z_t = \begin{bmatrix} 
Q_t^{\smash{1\!/2}} & 0_{n,m} \end{bmatrix}^{\trs} x_t + \begin{bmatrix} 0_{m,n} & R_t^{\smash{1\!/2}} \end{bmatrix}^{\trs} u_t, 
%\quad t\in\mathbb{N}_0,
\end{align}
\mc{where $A_t \in \mathbb{R}^{n \times n}$, $B_t \in \mathbb{R}^{n \times m}$, } $Q_t\in\mathbb{S}_+^n$, and $R_t\in\mathbb{S}_{\smash{++}}^m$. 

\begin{assumption} \label{asm:boundeddata}
The model data $A_t,B_t,Q_t,R_t$ \mc{and the inverses} $A_t^{-1}, R_t^{-1}$ are uniformly bounded across $t\in\mathbb{N}_0$.
\end{assumption}

\mc{Note that $A_t$ is non-singular whenever the model arises from discretization. Here it enables direct access to key  results on Riccati operator contraction~\cite{bougerol1993kalman}. Generalizing to singular $A_t$ (e.g., see~\cite{moore1980coping}) is beyond the current scope.}
\begin{assumption} \label{asm:uni_obs_ctr}
    The system model in~\eqref{eq:ltv_sys_w}--\eqref{eq:output_z_def}
    %system $(Q_t^{\smash{1\!/2}},A_t,B_t)_{t\in\mathbb{N}_0}$
    is uniformly observable and controllable in the following sense~\cite{anderson1981detectability}:
    \begin{align*}
    &(\exists d\in\mathbb{N})~ (\exists c\in\mathbb{R}_{>0})
    ~
    (\forall t \in \mathbb{N}_0) \\
    &~~~\left( \sum_{s=t}^{t+d-1} 
    \varPhi_{s,t}^{\trs} Q_s \varPhi_{s,t} \right) \succeq c\,I_n \preceq \left( \sum_{s=t}^{t+d-1} \varPhi_{s,t} B_sB_s^\trs \varPhi_{s,t}^\trs \right) ,  
    \end{align*}
    where 
    %$\mathbb{N}\!:=\!\{\!1,2,\ldots\!\}$, 
    $\varPhi_{t,t}\!:=\!I_n$, and
    $\varPhi_{s,t}\!:=\!A_{s-1}\cdot \ldots \cdot A_t$ for $s>t$.
    %denotes the state transition matrix.
    %, and $(\cdot)^{\trs}$ denotes~transpose.
\end{assumption}

While \mc{infinite-preview} state feedback controllers \mc{with} energy-gain \mc{guarantees} exist under uniform stabilizability and detectability, Assumption~\ref{asm:uni_obs_ctr} \mc{plays a role in the Riccati contraction based} synthesis \mc{of finite receding-horizon approximations, as elaborated in Sections~\ref{sec:mainresult} and~\ref{subsec:strcontract}}. 

The object of energy-gain control is to find a policy for $u=(u_t)_{t\in\mathbb{N}_0}$ \mc{such that the resulting map} from the disturbance input $w=(w_t)_{t\in\mathbb{N}_0}$
to the performance output $z=(z_t)_{t\in\mathbb{N}_0}$
is input-output stable over \mc{the space of finite energy signals}
\[\ell_2:=\{w=(w_t)_{t\in\mathbb{N}_0}~|~\|w\|_2^2:={\textstyle \sum_{t\in\mathbb{N}_0}} w_t^{\trs} w_t < +\infty\},\]
with \mc{prescribed worst-case} gain bound
$\gamma\in\mathbb{R}_{>0}$ \mc{in the sense that}
$(\forall w\in \ell_2)~J_\gamma(u,w) \leq 0$,
\mc{where} 
\begin{align}\label{eq:J_def}
J_\gamma(u,w) := \|z\|_2^2 - \gamma^2 \|w\|_2^2
= {\textstyle \sum_{t\in\mathbb{N}_0}}~ z_t^{\trs} z_t - \gamma^2 w_t^{\trs} w_t,
\end{align}
\mc{with} $z_t$ \mc{as per}~\eqref{eq:output_z_def} \mc{and}~\eqref{eq:ltv_sys_w} given $x_0=0$.
\mcr{This performance requirement amounts to} \mc{$\sup_{0\neq w\in\ell_2}\|z\|_2/\|w\|_2 \leq \gamma$,} \mcr{whereby} internal stability is implied under Assumptions~\ref{asm:boundeddata} and~\ref{asm:uni_obs_ctr}. 

The following result is
the standard \mc{infinite-horizon formulation} of a \mc{(strictly causal)} state feedback $\ell_2$ gain control policy; \mc{e.g., see}~\cite{green1994linear,hassibi1999indefinite},\mcr{\cite{katayama1994h}},~\cite[Sec.~2.4.2]{jintaothesis}.
\begin{theorem}\label{theorem:strictly_causal_gamma}
Given $\gamma \in \mathbb{R}_{>0}$, suppose the sequence  $(P_t)_{t\in\mathbb{N}_0}\subset\mathbb{S}_+^n$ \mc{satisfies the following}:
\begin{align}\label{eq:sc_P_ub}
& (\exists \varepsilon\in\mathbb{R}_{>0})~(\forall t \in  \mathbb{N}_0) ~~P_{t+1} - \gamma^2 I_n \preceq -\varepsilon I_n; ~~\text{ and } \\
\label{eq:RiccatiRecursion}
& (\forall t \in  \mathbb{N}_0) ~~ P_t = \mathcal{R}_{\gamma,t}(P_{t+1}),
\end{align}
where \mc{the $\gamma$-dependent time-varying $\ell_2$ gain Riccati operator for the system~\eqref{eq:ltv_sys_w}--\eqref{eq:output_z_def} is defined by}
\begin{align}\label{eq:RiccatiGamma} 
%\mc{P \mapsto~}
\mathcal{R}_{\gamma,t}(P) 
&:= Q_t\!+\!A_t^{\trs}PA_t\!-\!L_t(P)^{\trs}(M_{\gamma,t}(P))^{-1}L_t(P),
\end{align} 
with
\begin{align}
\begin{split}
\label{eq:LM_def}
L_t(P) &:= \begin{bmatrix}B_t & I_n\end{bmatrix}^{\trs} P A_t,\\
M_{\gamma,t}(P) &:= \begin{bmatrix}
R_t \!+\! B_t^{\trs}P B_t & B_t^{\trs}P \\
P B_t & P \!-\! \gamma^2I_n
\end{bmatrix}.
\end{split}
\end{align}
\mc{Further, in~\eqref{eq:ltv_sys_w}--\eqref{eq:output_z_def}, let}
\begin{align}\label{eq:h_inf_causal_policy}
u_t=\mathsf{u}_{\gamma,t}^{\mathsf{inf}}(x_t) := -K_{\gamma,t}(P_{t+1})\, x_t,
\quad t \in \mathbb{N}_0,
\end{align}
\mc{where}
\begin{align} 
\begin{split}
\label{eq:Kgam_nabla}
K_{\gamma,t}(P) &:= (\nabla_{\!\gamma,t}(P))^{-1} B_t^{\trs}
(P + P(\gamma^2I_n -  P)^{-1}P)
A_t,
\\
\nabla_{\!\gamma,t}(P) &:= R_t + B_t^{\trs}(P + P (\gamma^2 I_n-P)^{-1}P) B_t.
\end{split}
\end{align}
\mc{Then, with $J_\gamma$ as per~\eqref{eq:J_def},} $(\forall w \in \ell_2)~ J_\gamma(u,w) \leq 0$.
\end{theorem}

Related infinite-dimensional operator \mc{based} formulations of feedback controllers with \mc{guaranteed} $\ell_2$ gain  
are given in~\cite{halanay1994time,peters1997minimum,dullerud1999new}. As above, these also depend on all forward problem data. In particular, the state feedback policy in~\eqref{eq:h_inf_causal_policy} depends on $(Q_s,A_s,B_s,R_s)_{s\geq t}$ via the recursion~\eqref{eq:RiccatiRecursion}. The superscript $\mathsf{inf}$ emphasizes this \emph{infinite preview} dependence.

\mc{Infinite-preview dependence of the policy~\eqref{eq:h_inf_causal_policy} on the model parameters detracts from its practical applicability.}
Determining the \mc{hypothesized} solution $(P_t)_{t\in\mathbb{N}_0}$ of~\eqref{eq:RiccatiRecursion} is a challenge \mc{unless there is known structure, such as periodic  invariance~\cite{hench1994numerical,bittanti2009periodic}.} \mc{This motivates consideration of} state feedback $\ell_2$ gain control policy synthesis \mc{subject to finite preview of the model data in a receding-horizon fashion}. \mc{The proposed approach} is based on approximating $P_{t+1}$ in~\eqref{eq:h_inf_causal_policy}, as outlined \mc{in Section~\ref{subsec:maincontrib}}. \mc{It is established that} the error can be made arbitrarily small \mc{by using a sufficient number of model preview steps in the construction of the approximation. As such,} the development yields a practical method for computing the \mc{hypothesized} solution of~\eqref{eq:RiccatiRecursion} \mc{to desired accuracy at each time $t\in\mathbb{N}_0$. Continuity of closed-loop $\ell_2$ gain bound with respect to the approximation error then leads to the main receding-horizon controller synthesis result.}

\subsection{Main contributions} \label{subsec:maincontrib}

The main contribution is a finite receding-horizon preview synthesis of a state feedback controller 
%\mc{that approximates to} specified performance loss \mc{tolerance} 
that approximates the infinite-preview control policy\mcr{~\eqref{eq:h_inf_causal_policy}} in Theorem~\ref{theorem:strictly_causal_gamma}. The \mcr{prescribed} baseline 
%\mcr{disturbance energy-gain performance} 
%closed-loop 
performance bound $\gamma\in\mathbb{R}_{>0}$ 
%for the control policy
for the latter
is taken to be large enough for~\eqref{eq:RiccatiRecursion} to imply
$(P_{t})_{t\in\mathbb{N}_0}\subset\mathbb{S}_{\smash{++}}^n$, as elaborated in Remarks~\ref{rem:baseline} and~\ref{rem:PDP} \mc{in Sections~\ref{sec:mainresult} and~\ref{sec:proof}, respectively.}
 
As detailed in Section~\ref{sec:mainresult}, \mcr{given} \mc{performance loss} \mcr{tolerance}
$\beta \in \mathbb{R}_{>0}$, the proposed \mc{finite-preview} synthesis involves \mcr{approximating} $P_{t+1}$ in~\eqref{eq:h_inf_causal_policy}
by a positive definite $X_{t+1} \mc{\prec (\gamma+\beta)^2 I_n}$. 
%\mc{
The approximation is 
constructed by composing finitely many strictly contractive Riccati operators arising from a $d$-step lifting of~\eqref{eq:RiccatiRecursion}, in alignment with Assumption~\ref{asm:uni_obs_ctr}.
In this way, dependence on \mcr{the model parameters} 
%in~\eqref{eq:ltv_sys_w}--\eqref{eq:output_z_def} 
\mc{is confined to a} finite horizon ahead of $t\in\mathbb{N}_0$. \mc{The corresponding (strictly causal) state feedback control policy is} 
\begin{align}\label{eq:h_inf_approx_policy} u_t&=\mathsf{u}_{\smash{\gamma+\beta,t}}^{\smash{\mathsf{fin}}}(x_t) := -K_{\gamma+\beta,t}(X_{t+1})\, x_t, \quad t\in\mathbb{N}_0,
 \end{align}
where
$K_{\mc{\alpha},t}(X)= (\nabla_{\mc{\!\alpha},t}(X))^{-1}
B_t^{\trs}\left(X^{-1} - \mc{\alpha}^{-2}I_n\right)^{-1}\! A_t$, 
and $\nabla_{\!\mc{\alpha},t}(X)=R_t + B_t^\trs (X^{-1} - \mc{\alpha}^{-2}I_n)^{-1}B_t$, as per~\eqref{eq:Kgam_nabla} \mc{given that}
$X-X(X-\mc{\alpha}^2 I_n)^{-1}X=(X^{-1}-\mc{\alpha}^{-2}I_n)^{-1}$ for non-singular $X\prec \alpha^2 I_n$ by the Woodbury formula. The superscript $\mathsf{fin}$ emphasizes {\em finite preview} dependence on the model data.

The main result is formulated as Theorem~\ref{theorem:h_inf_T_lb}, in Section~\ref{sec:mainresult}. It %\mcr{gives} a sufficient number of 
gives a sufficient number of
preview steps in the proposed construction of each element of $(X_{t+1})_{t\in\mathbb{N}_0}$ for \mc{the resulting} \mcr{policy} $u=(u_t)_{t\in\mathbb{N}_0}=(\mathsf{u}_{\smash{\gamma+\beta,t}}^{\smash{\mathsf{fin}}}(x_t))_{t\in\mathbb{N}_0}$ \mcr{in~\eqref{eq:h_inf_approx_policy}} \mc{to achieve} $(\forall w \in \ell_2)~ J_{\gamma+\beta}(u,w) \leq 0$, with $J_{\gamma+\beta}$ as per~\eqref{eq:ltv_sys_w}--\eqref{eq:J_def}; \mcr{i.e., $\beta$ bounded energy-gain performance loss, relative to the infinite-preview policy~\eqref{eq:h_inf_causal_policy}.} The \mc{structured} proof developed in Section~\ref{sec:proof} combines Theorem~\ref{theorem:ric_gamma_contraction} on strict contraction of lifted Riccati operators, %\mc{composed to form each} $X_{t+1}$, 
and Theorem~\ref{theorem:delta_requirement} on  performance continuity. \mc{All other results, presented as lemmas, serve to establish these three main contributions.}

\subsection{Related work} \label{subsec:related}

A moving-horizon controller, \mc{for which} an infinite-horizon energy-gain \mc{bound exists}, is presented in~\cite{lall1995robust} for time-varying linear continuous-time systems. In the \mc{complementary setting of} discrete-time \mc{systems}, the \mc{distinctive} Riccati contraction based receding-horizon controller synthesis \mc{proposed here limits} performance degradation relative to the \mc{infinite-preview policy~\eqref{eq:h_inf_causal_policy}, thereby ensuring a prescribed worst-case energy-gain bound}. \mc{The receding-horizon controller in~\cite{essick2014control} also satisfies an $\ell_2$ gain specification for a class of switched systems. The finite structure underlying parameter switching plays a key role in refining the results of~\cite{dullerud1999new} to this class of time-varying systems. By contrast, the subsequent developments do not rely on such structural assumptions.}

In~\cite{goel2021regret,goel2023regret}, the setting is linear time-varying and discrete time, but the energy-gain performance horizon is finite. The focus is on performance degradation relative to non-causal preview of the disturbance. This so-called regret perspective is \mc{different to} the setup here, where the preview constraint relates to the availability of the model data in \mc{strictly causal} control policy synthesis with no preview of the disturbance. \mc{In~\cite{chen2025two}, limited preview of the cost parameters is considered from a regret perspective in a two-player linear quadratic game, but again the overall problem horizon is finite.} 

In receding horizon approximations of infinite-horizon optimal control policies, a terminal state penalty \mc{is typically imposed} in the finite-horizon optimization problem \mc{solved at} each step; e.g., see~\cite{kwon1978on,keerthi1988optimal,rawlings2017model}. \mc{For} time-invariant nonlinear continuous-time \mc{systems and no} disturbance, it is established in~\cite{jadbabaie2005stability} that with zero terminal cost, there exists a finite prediction horizon \mc{for which} stability \mc{is guaranteed}. In~\cite{grune2008infinite}, an infinite-horizon performance bound is also quantified for zero terminal penalty and set prediction horizon, with the added complication of state and input constraints, but nevertheless \mc{time-invariant} model parameters \mc{and no disturbance}. 

A receding-horizon feedback control policy that achieves a given $\ell_2$ gain bound appears in~\cite{goulart2009control} for constrained time-invariant linear discrete-time systems. The online synthesis involves the optimization of feedback policies over a finite prediction horizon~\cite{goulart2006optimization}, with terminal ingredients constructed via the algebraic $\ell_2$ gain Riccati equation for the unconstrained problem. This presumes constant problem data in a way that cannot be extended directly to a time-varying context without infinite preview \mc{of the model}. The same applies to the terminal ingredients in the receding-horizon regret-optimal schemes of~\cite{karapetyan2022regret,martin2025guarantees} for time-invariant systems, where \mc{again} regret relates to full preview of the disturbance. 

In a stationary setting, a contraction analysis of the Riccati operator for discrete-time block-update risk-sensitive filtering is developed in~\cite{levy2016contraction}. Under controllability and observability, the Riccati operator is shown to be strictly contractive with respect to Thompson's part metric, for a range of risk-sensitivity parameter values. The block-update implementation of the risk-sensitive filter is related to the type of lifting employed \mc{here}, where the time-varying setting \mc{is more general}, the Riemannian metric is used for contraction analysis, and the context is $\ell_2$ gain control.

\subsection{Organization} \label{subsec:map}

The article is organized as follows. The main result, \mc{outlined in Section~\ref{subsec:maincontrib}}, is formulated as Theorem~\ref{theorem:h_inf_T_lb} in Section~\ref{sec:mainresult}. \mc{This section includes development of the underlying finite-horizon lifting used} to obtain a one-step controllable and observable model, \mc{which enables direct application of an existing result on Riccati operator contraction in the synthesis of the approximating control policy}. A \mc{structured proof of the main result,} based on strict contraction of lifted $\ell_2$ gain Riccati operators, is developed in Section~\ref{sec:proof}. A numerical example is \mc{considered} in Section~\ref{sec:numex} \mc{to explore aspects of the main result}, followed by some concluding remarks in Section~\ref{sec:conc}.  

\section{Main result} \label{sec:mainresult}

\mc{With regard to the proposed receding-horizon $\ell_2$ gain controller synthesis,
% outlined in Section~\ref{subsec:maincontrib}
lifting the system model} at \mc{each} time $t\in\mathbb{N}_0$ \mc{enables} finite-preview construction of \mc{a positive definite} matrix $X_{t+1}$ to approximate $P_{t+1}$ \mc{in~\eqref{eq:h_inf_causal_policy}. Under Assumption~\ref{asm:uni_obs_ctr}, $d$-step lifting yields a one-step controllable and observable representation of the problem, unlocking existing theory on the contraction properties of Riccati operators~\cite{bougerol1993kalman}. This theory informs the construction of a suitable $X_{t+1}$ from a finite-horizon preview of the model data. A sufficient number of preview steps is identified in the main result, which is formulated as Theorem~\ref{theorem:h_inf_T_lb}, in terms of the lifted representation of the model~\eqref{eq:ltv_sys_w}--\eqref{eq:output_z_def} given in Lemmas~\ref{lemma:lifted_x}~and~\ref{lemma:lifted_z} below.}  

\mc{Given any} $d\in\mathbb{N}$ fixed in accordance with Assumption~\ref{asm:uni_obs_ctr}, and reference \mc{time} $t \in \mathbb{N}_0$, \mc{for every} 
$k\in\mathbb{N}_0$ define 
\begin{align} \label{eq:lifted_x_ch5}
\lftd{x}_{k|t} := x_{t+dk} \in \mc{\mathbb{R}^n},
\end{align}
the $d$-step lifted control input
\begin{align} \label{eq:lifted_u_ch5}
\lftd{u}_{k|t} :=
(u_{t+dk} \ , 
%\ u_{t+dk+1}, 
\ \cdots \ , 
\ u_{t+d(k+1)-1}) 
\in\mc{(\mathbb{R}^m)^d\sim\mathbb{R}^{md}},
\end{align}
the $d$-step lifted disturbance input
\begin{align}\label{eq:lifted_w_ch5}
\lftd{w}_{k|t} := 
(w_{t+dk} \ , 
%\ w_{t+dk+1}, 
\ \cdots \ , \ w_{t+d(k+1)-1}) \in \mc{
(\mathbb{R}^n)^d
\sim
\mathbb{R}^{nd}},
\end{align}
and
the $d$-step lifted performance output
\begin{align} \label{eq:lifted_z_ch5}
\lftd{z}_{k|t}:=(z_{t+dk}\ , \ \cdots\ ,\ z_{t+d(k+1)-1}) \in \mc{
(\mathbb{R}^{n+m})^d\sim
\mathbb{R}^{(n+m)d}}.
\end{align} 
\mc{The number of lifting steps $d$ is fixed in subsequent analysis, and as such, for convenience, the dependence on $d$ is suppressed in the notation.} From~\eqref{eq:ltv_sys_w},
\begin{align}\label{eq:dynamics_shorthand_ch5}
\mc{\varLambda_{k|t}}
\begin{bmatrix}
x_{t+dk} \\ \vdots \\ x_{t+d(k+1)-1} \\ \lftd{x}_{k+1|t}
\end{bmatrix}
=
\begin{bmatrix}
\lftd{x}_{k|t} \\ 0_{nd,1}
\end{bmatrix}
+ 
\mc{\varXi_{k|t}}
\lftd{u}_{k|t}
+ 
\begin{bmatrix}
0_{n,nd} \\ I_{nd}
\end{bmatrix}
\lftd{w}_{k|t},
\end{align}
where
\begin{align}
\mc{\varLambda_{k|t}} &:= I_{n(d+1)} -
\begin{bmatrix}
0_{n,nd} & 0_{n,n} \\
\diag(A_{t+dk}, \cdots, A_{t+d(k+1)-1}) & 0_{nd,n}
\end{bmatrix}, \label{eq:Ahat_ch5} \\
\mc{\varXi_{k|t}} &:= \begin{bmatrix}
0_{n,md} \\ \diag(B_{t+dk}, \cdots, B_{t+d(k+1)-1})
\end{bmatrix}. \label{eq:Bhat_ch5}
\end{align}
Similarly, from~\eqref{eq:output_z_def}, 
\begin{align}\label{eq:output_z_lifted_def}
\lftd{z}_{k|t} 
&= \begin{bmatrix} \mc{\varGamma_{k|t}} \\ 0_{md,n(d+1)} \end{bmatrix}
\begin{bmatrix}
x_{t+dk} \\ \vdots \\ x_{t+d(k+1)-1} \\ \lftd{x}_{k+1|t}
\end{bmatrix}
+
\begin{bmatrix}
0_{n(d+1),md} \\ \lftd{R}_{\smash{k|t}}^{\smash{1\!/2}}
\end{bmatrix}
\lftd{u}_{k|t},
\end{align}
where
\begin{align} \label{eq:Cbreve}
\mc{\varGamma_{k|t}} &:=
\begin{bmatrix}
\diag(Q_{t+dk}^{\smash{1\!/2}},\cdots,Q_{t+d(k+1)-1}^{\smash{1\!/2}}) & 0_{nd,n}
\end{bmatrix}, \\
\label{eq:hatR}
\lftd{R}_{k|t} &:= \diag(R_{t+dk},\cdots,R_{t+d(k+1)-1}).
\end{align}
On noting that $\mc{\varLambda_{k|t}}$ in~\eqref{eq:Ahat_ch5} is non-singular for all $t,k\in\mathbb{N}_0$, the \mc{next two lemmas are} a direct consequence of~\eqref{eq:dynamics_shorthand_ch5}.
\begin{lemma} \label{lemma:lifted_x}
Given  $u = (u_t)_{t\in\mathbb{N}_0}$ and $w=(w_t)_{t\in\mathbb{N}_0}$ \mc{in~\eqref{eq:ltv_sys_w}, with reference to~\eqref{eq:lifted_x_ch5}--\eqref{eq:lifted_w_ch5},
for every}
$t,k\in\mathbb{N}_0$,
\begin{align}\label{eq:lifted_dynamic_ch5}
 \lftd{x}_{k+1|t}
 &= \lftd{A}_{k|t} \lftd{x}_{k|t} + \lftd{B}_{k|t} \lftd{u}_{k|t} + \lftd{F}_{k|t} \lftd{w}_{k|t}, %
 \quad k\in\mathbb{N}_0,
\end{align}
%with $\lftd{x}_{0|t}=x_t$, $\lftd{u}_{k|t}$, and $\lftd{w}_{k|t}$ as per~\eqref{eq:lifted_x_ch5}--\eqref{eq:lifted_w_ch5}, 
where
\begin{align}
\lftd{A}_{k|t} &:= \begin{bmatrix}
0_{n,nd} & I_n
\end{bmatrix}
\mc{\varLambda_{\smash{k|t}}^{-1}}
\begin{bmatrix}
I_n \\ 0_{nd,n}
\end{bmatrix} \in \mathbb{R}^{n \times n}, \label{eq:phi_compact_ch5} \\
\lftd{B}_{k|t} &:= \begin{bmatrix}
0_{n,nd} & I_n
\end{bmatrix}
\mc{\varLambda_{\smash{k|t}}^{-1}} \mc{\varXi_{k|t}} \in \mathbb{R}^{n \times md}, \label{eq:gamma_compact_ch5} \\
\lftd{F}_{k|t} &:=
\begin{bmatrix}
0_{n,nd} & I_n
\end{bmatrix}
\mc{\varLambda_{\smash{k|t}}^{-1}}
\begin{bmatrix}
0_{n,nd} \\ I_{nd}
\end{bmatrix} \in \mathbb{R}^{n \times nd}, \label{eq:pi_compact_ch5}
\end{align}
$\mc{\varLambda_{k|t}}$ is defined in~\eqref{eq:Ahat_ch5}, and $\mc{\varXi_{k|t}}$ in~\eqref{eq:Bhat_ch5}.
\end{lemma}

In~\eqref{eq:gamma_compact_ch5}, $\lftd{B}_{k|t}$ is the $d$-step controllability matrix associated with the model data $(A_s, B_s)_{s \in \{t+dk, \ldots, t+d(k+1)-1\}}$.
Under Assumption~\ref{asm:uni_obs_ctr}, 
$(\exists c \in \mathbb{R}_{>0})~ (\forall t,k \in \mathbb{N}_0) ~ c I_n \preceq \lftd{B}_{k|t}\lftd{B}_{\smash{k|t}}^{\trs}$, making the lifted model~\eqref{eq:lifted_dynamic_ch5} one-step controllable. Further,
\begin{align*}
\lftd{F}_{k|t} \!=\! \begin{bmatrix}
\varPhi_{t+d(k+1),t+dk+1} & \cdots & \varPhi_{t+d(k+1),t+d(k+1)-1} & I_n
\end{bmatrix}\!,
\end{align*}
and thus, $(\exists \underline{c},\overline{c} \in \mathbb{R}_{>0}) ~ (\forall t,k \in \mathbb{N}_0) ~ \underline{c} I_n \preceq \lftd{F}_{k|t}\lftd{F}_{\smash{k|t}}^{\trs} \preceq \overline{c} I_n$. Moreover, 
$\lftd{A}_{k|t}=A_{t+d(k+1)-1} ~ \cdots ~ A_{t+dk} = \varPhi_{t+d(k+1),t+dk}$, $\lftd{A}_{\smash{k|t}}^{-1}$,
and $\lftd{B}_{k|t}$, are all uniformly bounded by Assumption~\ref{asm:boundeddata}.

\begin{lemma} \label{lemma:lifted_z}
Given $u = (u_t)_{t\in\mathbb{N}_0}$ and $w=(w_t)_{t\in\mathbb{N}_0}$ \mc{in~\eqref{eq:ltv_sys_w}, with reference to~\eqref{eq:lifted_x_ch5}--\eqref{eq:lifted_z_ch5},
for every} $t,k\in\mathbb{N}_0$, 
\begin{align}\label{eq:output_z_lifted}
\lftd{z}_{k|t} = \begin{bmatrix} \lftd{C}_{k|t} \\ 0_{md,n(d+1)} \end{bmatrix} \! \lftd{x}_{k|t}
+ \begin{bmatrix} \lftd{D}_{k|t} \\ \lftd{R}_{\smash{k|t}}^{\smash{1\!/2}} \end{bmatrix} \! \lftd{u}_{k|t}
+ \begin{bmatrix} \lftd{E}_{k|t} \\ 0_{md,n(d+1)} \end{bmatrix} \! \lftd{w}_{k|t},
\end{align}
% with $\lftd{u}_{k|t}$, $\lftd{w}_{k|t}$, and $\lftd{z}_{k|t}$ as per~\eqref{eq:lifted_u_ch5}--\eqref{eq:lifted_z_ch5} \mc{and~\eqref{eq:output_z_lifted_def}}, and $\lftd{x}_{\mc{k|t}}$ as per~\eqref{eq:lifted_dynamic_ch5},
where
\begin{align}
\lftd{C}_{k|t} &:= \mc{\varGamma_{k|t}} \mc{\varLambda_{\smash{k|t}}}^{-1} \begin{bmatrix}
I_n \\ 0_{nd,n}
\end{bmatrix} \in \mathbb{R}^{nd \times n}, \label{eq:C_ch5} \\
\lftd{D}_{k|t} &:= \mc{\varGamma_{k|t}} \mc{\varLambda_{\smash{k|t}}^{-1}} \mc{\varXi_{k|t}} \in \mathbb{R}^{nd \times md}, \label{eq:D_ch5} \\
\lftd{E}_{k|t} &:= \mc{\varGamma_{k|t}} \mc{\varLambda_{\smash{k|t}}^{-1}}
\begin{bmatrix} 0_{n,nd} \\ I_{nd} \end{bmatrix} \in \mathbb{R}^{nd \times nd}, \label{eq:E_ch5}
\end{align}
$\mc{\varLambda_{k|t}}$ is defined in~\eqref{eq:Ahat_ch5}, $\mc{\varXi_{k|t}}$ in~\eqref{eq:Bhat_ch5}, $\mc{\varGamma_{k|t}}$ in~\eqref{eq:Cbreve}, and $\lftd{R}_{k|t}$ in~\eqref{eq:hatR}.
\end{lemma}

In~\eqref{eq:C_ch5}, $\lftd{C}_{k|t}$ is the $d$-step observability matrix associated with the model data $(Q_s^{\smash{1\!/2}}, A_s)_{s \in \{t+dk, \ldots, t+d(k+1)-1\}}$.
Under Assumption~\ref{asm:uni_obs_ctr}, 
$(\exists c\in\mathbb{R}_{>0})~(\forall t,k \in \mathbb{N}_0)~c I_n \preceq \lftd{C}_{\smash{k|t}}^{\trs}\lftd{C}_{k|t}$ making the lifted model \eqref{eq:lifted_dynamic_ch5} with output~\eqref{eq:output_z_lifted} one-step observable.
Moreover, $\lftd{C}_{k|t}$, $\lftd{D}_{k|t}$, and $\lftd{E}_{k|t}$ are all uniformly bounded by Assumption~\ref{asm:boundeddata}.

\mc{The following finite-preview construction of \mc{a matrix} $X_{t+1}$ to approximate $P_{t+1}$ in~\eqref{eq:h_inf_causal_policy} involves a final transformation of the model data. This arises to accommodate the direct dependence of $\lftd{z}_{k|t}$ on $\lftd{w}_{k|t}$ in the supporting analysis elaborated in Section~\ref{subsec:liftricc}.} For $t,k\in\mathbb{N}_0$, with $\lftd{A}_{k|t},\lftd{B}_{k|t}, \lftd{C}_{k|t},\lftd{F}_{k|t},\lftd{D}_{k|t},\lftd{E}_{k|t}$, and $\lftd{R}_{k|t}$ \mc{as given} in Lemmas~\ref{lemma:lifted_x} and~\ref{lemma:lifted_z}, define
\begin{align}
\tilde{\lftd{B}}_{k|t} &:= \begin{bmatrix}
\lftd{B}_{k|t} & \lftd{F}_{k|t}
\end{bmatrix}, \label{eq:btilch5} \\
\label{eq:rtilch5} 
\tilde{\lftd{R}}_{k|t} &:=
\begin{bmatrix}
\lftd{R}_{k|t} + \lftd{D}_{\smash{k|t}}^{\trs} \lftd{D}_{k|t} & \lftd{D}_{\smash{k|t}}^{\trs} \lftd{E}_{k|t} \\
\lftd{E}_{\smash{k|t}}^{\trs} \lftd{D}_{k|t} & \lftd{E}_{\smash{k|t}}^{\trs} \lftd{E}_{k|t} - \gamma^2 I_{nd}
\end{bmatrix},
\end{align}
\mc{with $\gamma\in\mathbb{R}_{>0}$} such that $\tilde{\lftd{R}}_{\mc{k|t}}$ is non-singular, 
\begin{align}
    \tilde{\lftd{Q}}_{k|t} &:= \lftd{C}_{\smash{k|t}}^{\trs}\lftd{C}_{k|t} - \lftd{C}_{\smash{k|t}}^{\trs} \begin{bmatrix} \lftd{D}_{k|t} & \lftd{E}_{k|t} \end{bmatrix}
\tilde{\lftd{R}}_{\smash{k|t}}^{-1}
\begin{bmatrix} \lftd{D}_{\smash{k|t}}^{\trs} \\ \lftd{E}_{\smash{k|t}}^{\trs} \end{bmatrix} \lftd{C}_{k|t},   
\label{eq:qtilch5} 
\\
\tilde{\lftd{A}}_{k|t} &:= \lftd{A}_{k|t} - \begin{bmatrix} \lftd{B}_{k|t} & \lftd{F}_{k|t} \end{bmatrix} 
\tilde{\lftd{R}}_{\smash{k|t}}^{-1}
\begin{bmatrix} \lftd{D}_{\smash{k|t}}^{\trs} \\ \lftd{E}_{\smash{k|t}}^{\trs} \end{bmatrix} \lftd{C}_{k|t}. \label{eq:atilch5}
\end{align} 
The dependence on $\gamma$ \mc{is} suppressed for convenience. \mc{Note that non-singularity of $\tilde{R}_{k|t}$ implies non-singularity of $\tilde{A}_{k|t}$ by Lemma~\ref{lemma:lifted_A_nonsingular} in the Appendix. Finally,} given $T\in\mathbb{N}$,
let
\begin{align}\label{eq:h_inf_terminal}
\tilde{X}_{t+1} := \tilde{\lftd{Q}}_{T|t+1} + \tilde{\lftd{A}}_{\smash{T|t+1}}^{\trs} (\tilde{\lftd{B}}_{T|t+1} \tilde{\lftd{R}}_{\smash{T|t+1}}^{-1} \tilde{\lftd{B}}_{\smash{T|t+1}}^{\trs})^{-1} \tilde{\lftd{A}}_{T|t+1},
\end{align}
subject to \mc{non-singularity of} $\tilde{\lftd{B}}_{T|t+1} \tilde{\lftd{R}}_{\smash{T|t+1}}^{-1} \tilde{\lftd{B}}_{\smash{T|t+1}}^{\trs}$. If, \mc{in addition},  $\tilde{\lftd{B}}_{T|t+1} \tilde{\lftd{R}}_{\smash{T|t+1}}^{-1} \tilde{\lftd{B}}_{\smash{T|t+1}}^{\trs}\in\mathbb{S}_{\smash{++}}^n$ and $\tilde{\lftd{Q}}_{T|t+1}\in\mathbb{S}_{\smash{++}}^n$, \mc{which would be infeasible were the lifted model not one-step controllable and observable}, then $\tilde{X}_{t+1}\in\mathbb{S}_{\smash{++}}^n$. \mc{Given} this, $(X_{t+1})_{t\in\mathbb{N}_0}\subset\mathbb{S}_{\smash{++}}^n$ when defined according to
\begin{align}\label{eq:h_inf_approx_X}
X_{t+1} 
&:= 
\tilde{\lftd{\mathcal{R}}}_{\gamma,0|t+1}
\circ \cdots \circ 
\tilde{\lftd{\mathcal{R}}}_{\gamma,T-1|t+1}(\tilde{X}_{t+1}),
\end{align}
%\mc{with}
%$\tilde{X}_{t+1}$ \mc{as in}~\eqref{eq:h_inf_terminal}, 
\mc{where the lifted} 
Riccati operators \mc{are defined as follows:}
\begin{align}\label{eq:lifted_ric_0}
&\tilde{\lftd{\mathcal{R}}}_{\gamma,k|t}(X) \nonumber \\
&\!:= 
\tilde{\lftd{Q}}_{k|t} \!+\! \tilde{\lftd{A}}_{\smash{k|t}}^{\trs}( X \!-\! X\tilde{\lftd{B}}_{k|t}(\tilde{\lftd{R}}_{k|t}\!+\! \tilde{\lftd{B}}_{\smash{k|t}}^{\trs}X\tilde{\lftd{B}}_{k|t})^{-1} \tilde{\lftd{B}}_{\smash{k|t}}^{\trs}X)\tilde{\lftd{A}}_{k|t}.
\end{align} 
\mc{This} construction of $X_{t+1}$ involves $d\cdot(T+1)$ \mc{steps of model data} ahead of $t\in\mathbb{N}_0$. The main result, \mcr{formulated} below, establishes a sufficient number of \mcr{lifted} preview steps \mcr{$T$} for the corresponding control policy~\eqref{eq:h_inf_approx_policy} to meet a specified performance loss bound $\beta\in\mathbb{R}_{>0}$, relative to the infinite-preview policy~\eqref{eq:h_inf_causal_policy} \mcr{with prescribed} baseline \mcr{gain} bound $\gamma$. 
\begin{theorem}\label{theorem:h_inf_T_lb}
Given $\gamma\in\mathbb{R}_{>0}$, suppose:
\begin{enumerate}
    \item[1)] \mc{the hypothesis in Theorem~\ref{theorem:strictly_causal_gamma} holds};~and
    \item[2)] there exist $\underline{c},\overline{c}\in\mathbb{R}_{>0}$, such that for all $t\in\mathbb{N}_0$, 
    \begin{enumerate}
        \item[a)] 
$\underline{c}I_{(m+n)d} \preceq \tilde{\lftd{R}}_{\smash{0|t}}^\trs \tilde{\lftd{R}}_{0|t} \preceq \overline{c}I_{(m+n)d}$, 
\item[b)] $\underline{c}I_n \preceq \tilde{\lftd{B}}_{0|t}\tilde{\lftd{R}}_{\smash{0|t}}^{-1}\tilde{\lftd{B}}_{\smash{0|t}}^{\trs} \preceq \overline{c}I_n$, and
\item[c)] $\underline{c}I_n \preceq \tilde{\lftd{Q}}_{0|t} \preceq \overline{c}I_n$,
\end{enumerate} 
with $\tilde{\lftd{B}}_{0|t}$ as per~\eqref{eq:btilch5}, and $\gamma$-dependent $\tilde{\lftd{R}}_{0|t}$ and $\tilde{\lftd{Q}}_{0|t}$ as per~\eqref{eq:rtilch5} and~\eqref{eq:qtilch5}, respectively.
\end{enumerate} 
Further, \mc{given $\beta\in\mathbb{R}_{>0}$}, let $u=(u_t)_{t\in\mathbb{N}_0}\!=\!(\mathsf{u}_{\smash{\gamma+\beta,t}}^{\smash{\mathsf{fin}}}(x_t))_{t\in\mathbb{N}_0}$ \mc{as per} the finite-preview state feedback control policy~\eqref{eq:h_inf_approx_policy} for the system~\eqref{eq:ltv_sys_w}--\eqref{eq:output_z_def}, with $(X_{t+1})_{t\in\mathbb{N}_0}$ \mc{constructed according to}~\eqref{eq:h_inf_approx_X} \mc{for any selection of} 
\begin{align}\label{eq:preview_lb}
T > \underline{T} := \log \Big( ~\log \big( ~(\underline{\alpha} 
 + 1)^{1/\overline{\delta}}~
\big)~\Big)  \Big/ \log(\overline{\rho}),
\end{align}
where 
$\underline{\alpha}:=(\gamma^{-2}-(\gamma+\beta)^{-2})\cdot \underline{\kappa}$,
\begin{align}
\label{eq:kappalo}
\underline{\kappa}&:= \inf_{t\in\mathbb{N}_0} \lambda_{\min}(\tilde{\lftd{Q}}_{0|t}) >0, \\
\overline{\delta} &:= \!
\sqrt{n} \log \!\left(
\sup_{t\in\mathbb{N}_0} \!
\frac{\lambda_{\max}(\tilde{\lftd{Q}}_{0|t} \!+\! \tilde{\lftd{A}}_{\smash{0|t}}^{\trs} (\tilde{\lftd{B}}_{0|t} \tilde{\lftd{R}}_{\smash{0|t}}^{-1} \tilde{\lftd{B}}_{\smash{0|t}}^{\trs})^{-1}\!\tilde{\lftd{A}}_{0|t})}{\lambda_{\min}(\tilde{\lftd{Q}}_{0|t})} \right) 
\nonumber \\ 
&\, < +\infty, \label{eq:deltaup}
\\
\label{eq:sup_rho_ch5}
\overline{\rho} &:= \sup_{t \in \mathbb{N}_0} 1/( 1 + \tilde{\omega}_{t} ) <1, \\
\label{eq:omtil}
\tilde{\omega}_{t} &:= \!\frac{\lambda_{\min}( \tilde{\lftd{Q}}_{0|t} \!+\! \tilde{\lftd{Q}}_{0|t}\tilde{\lftd{A}}_{\smash{0|t}}^{-1} \tilde{\lftd{B}}_{0|t}\tilde{\lftd{R}}_{\smash{0|t}}^{-1}\tilde{\lftd{B}}_{\smash{0|t}}^{\trs} (\tilde{\lftd{A}}_{\smash{0|t}}^{\trs})^{-1}\tilde{\lftd{Q}}_{0|t})}
{\lambda_{\max}(\tilde{\lftd{Q}}_{0|t} \!+\!  \tilde{\lftd{A}}_{\smash{0|t}}^{\trs}(\tilde{\lftd{B}}_{0|t}\tilde{\lftd{R}}_{\smash{0|t}}^{-1}\tilde{\lftd{B}}_{\smash{0|t}}^{\trs})^{-1} \tilde{\lftd{A}}_{0|t})},\!
\end{align}
\mc{and $\tilde{\lftd{A}}_{0|t}$ is as defined in~\eqref{eq:atilch5}.}
Then, with $J_{\gamma+\beta}$ as per~\eqref{eq:J_def},
$(\forall w\in\ell_2)~J_{\gamma+\beta}(u,w) \leq 0$.
\end{theorem}
\begin{remark} \label{rem:baseline}
Parts 1) and 2) of the hypothesis in Theorem~\ref{theorem:h_inf_T_lb} amount to \mc{considering a} sufficiently large gain bound $\gamma$ for the baseline infinite-preview control policy~\eqref{eq:h_inf_causal_policy}. 
Given~\eqref{eq:sc_P_ub}, the Schur complement of 
$M_{\gamma,t}(P_{t+1})$, \mc{as} defined in~\eqref{eq:LM_def}, is given by
$R_t \!+\! B_t^\prime (P_{t+1}\! -\! P_{t+1}(P_{t+1}\!-\! \gamma^2 I_n)^{-1} P_{t+1}) B_t \!=\! \nabla_{\!\gamma,t}(P_{t+1})$; cf.~\eqref{eq:Kgam_nabla}. Further,
$P_{t+1} \!+\! P_{t+1}(\gamma^2 I_n\!-\! P_{t+1})^{-1} P_{t+1}\in\mathbb{S}_{\smash{+}}^n$, whereby  
$\nabla_{\!\gamma,t}(P_{t+1})\in\mathbb{S}_{\smash{++}}^m$.
%, as  $R_t\in\mathbb{S}_{\smash{++}}^m$.
Thus, $M_{\gamma,t}(P_{t+1})$ is non-singular for $t\in\mathbb{N}_0$;~\mc{i.e.,} Riccati recursion~\eqref{eq:RiccatiRecursion} is well posed.
In fact, $(P_t)_{t\in\mathbb{N}_0}\subset\mathbb{S}_{\smash{++}}^n$ under part 2) of the hypothesis, \mc{as elaborated in} Remark~\ref{rem:PDP}, \mc{Section~\ref{subsec:liftricc}}. 
\mc{It is important to note that the} construction of $(X_{t+1})_{t\in\mathbb{N}_0}$ \mc{according to}~\eqref{eq:h_inf_approx_X} does not \mc{involve knowledge of} $(P_t)_{t\in\mathbb{N}_0}$. \mc{The dependence on $\gamma$ relates to the performance loss perspective used to assess the resulting finite-preview controller~\eqref{eq:h_inf_approx_policy}.}
%, relative 
%to~\eqref{eq:h_inf_causal_policy}.}
\end{remark}
\begin{remark} \label{rem:deadbeat}
\mc{Under} part 2) of the hypothesis 
in Theorem~\ref{theorem:h_inf_T_lb},  $\tilde{X}_{t+1}\in\mathbb{S}_{++}^n$ in~\eqref{eq:h_inf_terminal}, and $\mathcal{R}_{\gamma,k|t}(X)\in\mathbb{S}_{\smash{++}}^n$ over $X\in\mathbb{S}_{\smash{++}}^n$, 
as elaborated in Remark~\ref{rem:part2forallt}, \mc{Section~\ref{subsec:liftricc}. Thus,} $X_{t+1}\in\mathbb{S}_{\smash{++}}^n$ in~\eqref{eq:h_inf_approx_X}. \mc{Using one-step controllability and observability of the lifted model, it is established in~\cite[Sec.~5.4.2]{jintaothesis} that the three components of this part of the} hypothesis \mcr{are necessary conditions} when for every $t\in\mathbb{N}_0$, the lifted deadbeat policy
\begin{align}
\lftd{u}_{k|t}
\!=\! -\lftd{B}_{\smash{k|t}}^{\trs}(\lftd{B}_{k|t} \lftd{B}_{\smash{k|t}}^{\trs})^{-1}(\lftd{A}_{k|t} \lftd{x}_{k|t} \!+\! \lftd{F}_{k|t} \lftd{w}_{k|t}),~k\in\mathbb{N}_0,
\label{eq:deadbeat_ch5}
\end{align} 
%k\in\mathbb{N}_0$, 
achieves the \mcr{baseline} gain bound $\gamma$, \mc{given $u_s=0$ and $w_s=0$ for $s<t$}. This full-information feedback control policy is not causal in the original time domain since \mc{in addition to model data} it involves $d$-step preview of the disturbance.
\end{remark}
\begin{remark} \label{rem:ddep}
    The quantity $\underline{\kappa}$ in~\eqref{eq:kappalo} 
    depends on the fixed number of steps $d$ used to \mc{lift the model to the form}~\eqref{eq:lifted_dynamic_ch5} and~\eqref{eq:output_z_lifted}. \mc{This number} is not uniquely determined by complying with Assumption~\ref{asm:uni_obs_ctr}, \mc{as} any greater value also complies. Likewise, $\overline{\delta}$ in~\eqref{eq:deltaup}, which bounds the Riemannian distance between $P_{t+dT+1}$ and $\tilde{X}_{t+1}$, and $\overline{\rho}$ in~\eqref{eq:sup_rho_ch5}, which bounds the contraction rate of the lifted Riccati operator in~\eqref{eq:lifted_ric_0},
    both depend on $d$. \mc{These quantities also all depend on $\gamma$.}
\end{remark}
\begin{remark} \label{rem:bounds}
    In principle, the quantities~\eqref{eq:kappalo}--\eqref{eq:omtil} can be replaced by corresponding lower or upper estimates determined from uniform bounds on the problem data \mc{in accordance with} Assumption~\ref{asm:boundeddata}. For periodic systems, these quantities, and $(P_t)_{t\in\mathbb{N}_0}$ for that matter~\cite{hench1994numerical,bittanti2009periodic}, can be computed directly from the finite amount of model data. This facilitates numerical \mc{investigation} of the kind \mc{presented} in Section~\ref{sec:numex} to assess the conservativeness of Theorem~\ref{theorem:h_inf_T_lb}, which is only sufficient to guarantee the performance loss bound $\beta$, relative to a suitably large baseline gain bound~$\gamma$.
    \end{remark}

\section{Proof of the main result} \label{sec:proof}

This section encompasses a \mc{structured} proof of Theorem~\ref{theorem:h_inf_T_lb}. \mc{This main result is established by combining  Theorems~\ref{theorem:ric_gamma_contraction} and~\ref{theorem:delta_requirement}, which are formulated below.} 

In particular, the \mc{proof relies on} the {\em strict} contraction of the lifted Riccati operators composed in~\eqref{eq:h_inf_approx_X} \mc{to form $(X_{t+1})_{t\in\mathbb{N}_0}$}.
%, with respect to the Riemannian metric on $\mathbb{S}_{\smash{++}}^n$. 
This key property holds under part 2) of the hypothesis on the baseline \mc{gain bound} $\gamma\in\mathbb{R}_{>0}$ \mc{in Theorem~\ref{theorem:h_inf_T_lb}}, \mc{as summarized in} Theorem~\ref{theorem:ric_gamma_contraction}, \mc{which is} proved in Section~\ref{subsec:strcontract}; \mc{recall that non-singularity of $\gamma$-dependent $\tilde{\lftd{R}}_{k|t}$ in~\eqref{eq:rtilch5} implies non-singularity of $\tilde{\lftd{A}}_{k|t}$ in~\eqref{eq:atilch5}, by Lemma~\ref{lemma:lifted_A_nonsingular} in the Appendix.} 
% by building upon several lemmas about lifted Riccati operators presented in Section~\ref{subsec:liftricc}\mc{:}
\begin{theorem} \label{theorem:ric_gamma_contraction}
Given $t,k\in\mathbb{N}_0$, and $\gamma\in\mathbb{R}_{>0}$, suppose $\tilde{\lftd{R}}_{k|t}$ in~\eqref{eq:rtilch5} is
non-singular, \mc{and} with reference to~\eqref{eq:btilch5} \mc{and \eqref{eq:qtilch5}}, suppose $\tilde{\lftd{B}}_{k|t}\tilde{\lftd{R}}_{\smash{k|t}}^{-1}\tilde{\lftd{B}}_{\smash{k|t}}^{\trs} \in \mathbb{S}_{\smash{++}}^n$ and $\tilde{\lftd{Q}}_{k|t} \in \mathbb{S}_{\smash{++}}^n$. Then, for any 
$X,P \in \mathbb{S}_{\smash{++}}^n$,
\begin{align}\label{eq:rho_ch5}
\delta(\tilde{\lftd{\mathcal{R}}}_{\gamma,k|t}(X),\tilde{\lftd{\mathcal{R}}}_{\gamma,k|t}(P)) \leq \tilde{\rho}_{k|t} \cdot \delta(X,P)
\end{align}
with $\tilde{\rho}_{k|t} := \tilde{\zeta}_{k|t}\big/(\tilde{\zeta}_{k|t} + \tilde{\epsilon}_{k|t}) < 1$, where
$\tilde{\lftd{\mathcal{R}}}_{\gamma,k|t}(\cdot)$ is defined in~\eqref{eq:lifted_ric_0},
\begin{align}
\tilde{\zeta}_{k|t} &:= 
1\big/\lambda_{\min}( \tilde{\lftd{Q}}_{k|t} + \tilde{\lftd{Q}}_{k|t}\tilde{\lftd{A}}_{\smash{k|t}}^{-1} \tilde{\lftd{B}}_{k|t}\tilde{\lftd{R}}_{\smash{k|t}}^{-1}\tilde{\lftd{B}}_{\smash{k|t}}^{\trs} (\tilde{\lftd{A}}_{\smash{k|t}}^{\trs})^{-1}\tilde{\lftd{Q}}_{k|t}),
\label{eq:h_inf_zeta}\\
\tilde{\epsilon}_{k|t} &:= 1 \big/ 
\lambda_{\max}(\tilde{\lftd{Q}}_{k|t} +  \tilde{\lftd{A}}_{\smash{k|t}}^{\trs}(\tilde{\lftd{B}}_{k|t}\tilde{\lftd{R}}_{\smash{k|t}}^{-1}\tilde{\lftd{B}}_{\smash{k|t}}^{\trs})^{-1} \tilde{\lftd{A}}_{k|t}),
\label{eq:h_inf_epsilon}
\end{align}
and $\delta(\cdot,\cdot)$ is the Riemannian \mc{metric on $\mcr{\mathbb{S}_{\smash{++}}^n}$; see}~\eqref{eq:riemannian}.
\end{theorem}

\mc{The other} key ingredient 
%in the proof of Theorem~\ref{theorem:h_inf_T_lb} 
\mc{pertains} to continuity of the baseline \mc{gain bound associated with the} infinite-preview policy~\eqref{eq:h_inf_causal_policy}.  
\begin{theorem}\label{theorem:delta_requirement}
Given $\gamma\in\mathbb{R}_{>0}$, suppose 
$(P_t)_{t\in\mathbb{N}_0}\subset\mathbb{S}_{\smash{++}}^n$ satisfies~\eqref{eq:sc_P_ub} and~\eqref{eq:RiccatiRecursion} in addition~to
$\underline{\lambda} := \inf_{t \in \mathbb{N}_0} \lambda_{\min}(P_t)>0$.
Further, given 
%performance loss specification 
$\beta\in\mathbb{R}_{>0}$, and $(X_{t+1})_{t\in\mathbb{N}_0}\subset\mathbb{S}_{\smash{++}}^n$, suppose
\begin{align}\label{eq:X_requirement_final}
&(\exists \varepsilon\in\mathbb{R}_{>0})~
(\forall t \in \mathbb{N}_0) \nonumber \\
&
P_{t+1} \prec X_{t+1}  ~\land~ 
\delta(X_{t+1},P_{t+1}) \leq 
\log((\mc{\eta}-\varepsilon) \underline{\lambda} + 1),
\end{align}
where \mcr{$\land$ denotes logical conjunction},
\begin{align}
\label{eq:alpha_ch5}
\mc{\eta} := \gamma^{-2} - (\gamma+\beta)^{-2}>0,
\end{align}
and $\delta(\cdot,\cdot)$ is the Riemannian distance; see~\eqref{eq:riemannian}. \mc{Finally,} let $u=(u_t)_{t\in\mathbb{N}_0}\!=\!(\mathsf{u}_{\smash{\gamma+\beta,t}}^{\smash{\mathsf{fin}}}(x_t))_{t\in\mathbb{N}_0}$ \mc{as per} the finite-preview state feedback control policy~\eqref{eq:h_inf_approx_policy} for the system~\eqref{eq:ltv_sys_w}--\eqref{eq:output_z_def}.
Then, with $J_{\gamma+\beta}$ as per~\eqref{eq:J_def},  $(\forall w\in\ell_2)~J_{\gamma+\beta}(u,w) \leq 0$.
\end{theorem}

\mc{Relationships between the Riccati operator~\eqref{eq:lifted_ric_0} for the lifted model~\eqref{eq:lifted_dynamic_ch5} with~\eqref{eq:output_z_lifted}, and the $\ell_2$ gain Riccati operator~\eqref{eq:RiccatiGamma} for~\eqref{eq:ltv_sys_w}--\eqref{eq:output_z_def}, are established in Section~\ref{subsec:liftricc}. This forms the basis for the proof of Theorem~\ref{theorem:delta_requirement} in Section~\ref{subsec:approx}. The development of these relationships also highlights interesting links between Riccati recursions and Schur decompositions of quadratic forms in $\ell_2$ gain analysis. In Section~\ref{subsec:proof}, where Theorem~\ref{theorem:h_inf_T_lb} is proved, it is shown via Theorem~\ref{theorem:ric_gamma_contraction} that with sufficiently large $T\in\mathbb{N}$, the construction of $(X_{t+1})_{t\in\mathbb{N}_0}$ in~\eqref{eq:h_inf_approx_X} satisfies the hypothesis~\eqref{eq:X_requirement_final}.}

\subsection{\mc{Lifted Riccati operator contraction}}\label{subsec:strcontract}

Applying~\cite[Thm.~1.7]{bougerol1993kalman}, the hypothesis in Theorem~\ref{theorem:ric_gamma_contraction} is shown to be sufficient for {\em strict} contraction of \mc{the lifted Riccati operator $X\mapsto  \tilde{\lftd{\mathcal{R}}}_{\gamma,k|t}(X)$ given by~\eqref{eq:lifted_ric_0}}. \mc{Without Assumption~\ref{asm:uni_obs_ctr}, the part of the hypothesis regarding $\tilde{\lftd{B}}_{k|t}\tilde{\lftd{R}}_{\smash{k|t}}^{-1}\tilde{\lftd{B}}_{\smash{k|t}}^{\trs}$ and $\tilde{\lftd{Q}}_{k|t}$ is infeasible, and the contraction may be non-strict.} The Riemannian distance between $X,P\in\mathbb{S}_{\smash{++}}^n$ is defined \mc{as follows~\cite{bhatia2007positive}}:
\begin{align} \label{eq:riemannian}
  \delta(X,P):= {\textstyle \sqrt{\sum_{i\in\{1,\ldots,n\}} (\log(\lambda_i))^2}},   
\end{align} 
where $\{\lambda_1,\ldots,\lambda_n\}=\mathrm{spec}(XP^{-1})$ is the spectrum of $XP^{-1}$ (i.e., the collection of all eigenvalues, including multiplicity.) The latter coincides with $\mathrm{spec}(P^{-\smash{1\!/2}}XP^{-\smash{1\!/2}})\subset\mathbb{R}_{>0}$ because $\mathrm{spec}(YZ)\cup\{0\}=\mathrm{spec}(ZY)\cup\{0\}$ for all square $Y,Z$, and $\{0\}\cap \mathrm{spec}(XP^{-1})=\emptyset =\{0\}\cap\mathrm{spec}(P^{-\smash{1\!/2}}XP^{-\smash{1\!/2}})$. Indeed,    
$\delta(X,P)=\delta(X^{-1},P^{-1})=\delta(P^{-1},X^{-1})=\delta(P,X)$ as $\lambda_i\in\mathrm{spec}(XP^{-1})=\mathrm{spec}(P^{-1}X)$ implies $1/\lambda_i\in\mathrm{spec}(X^{-1}P)=\mathrm{spec}(PX^{-1})$, and $(\log(\lambda_i))^2=(\log(1/\lambda_i))^2$.

{\noindent\hspace{1em}{\bf\itshape Proof of Theorem~\ref{theorem:ric_gamma_contraction}:}
%With $\tilde{\lftd{R}}_{k|t}$ non-singular,
%$\tilde{B}_{k|t}\tilde{R}_{\smash{k|t}}^{-1}\tilde{B}_{\smash{k|t}}^{\trs} \in \mathbb{S}_{\smash{++}}^n$, and $\tilde{Q}_{k|t} \in \mathbb{S}_{\smash{++}}^n$, 
Under the stated hypothesis, it follows by Remark~\ref{rem:prePDP} (Section~\ref{subsec:liftricc}), and Lemma~\ref{lemma:lifted_A_nonsingular} (Appendix), that
\begin{align*}
\tilde{\lftd{\mathcal{R}}}_{\gamma,k|t}(P)
&= 
\tilde{\lftd{Q}}_{k|t} + \tilde{\lftd{A}}_{\smash{k|t}}^{\trs}P(\tilde{\lftd{A}}_{\smash{k|t}}^{-1} + \tilde{\lftd{A}}_{\smash{k|t}}^{-1}\tilde{\lftd{B}}_{k|t}\tilde{\lftd{R}}_{\smash{k|t}}^{-1}\tilde{\lftd{B}}_{\smash{k|t}}^{\trs} P)^{-1}\\
&= \left(\! \tilde{\lftd{Q}}_{k|t}(\tilde{\lftd{A}}_{\smash{k|t}}^{-1} \!+\! \tilde{\lftd{A}}_{\smash{k|t}}^{-1}\tilde{\lftd{B}}_{k|t} \tilde{\lftd{R}}_{\smash{k|t}}^{-1} \tilde{\lftd{B}}_{\smash{k|t}}^{\trs} P)\! +\! \tilde{\lftd{A}}_{\smash{k|t}}^{\trs}P\! \right) \\
&\qquad\qquad\qquad \times \left(\!\tilde{\lftd{A}}_{\smash{k|t}}^{-1} \!+\! \tilde{\lftd{A}}_{\smash{k|t}}^{-1}\tilde{\lftd{B}}_{k|t} \tilde{\lftd{R}}_{\smash{k|t}}^{-1} \tilde{\lftd{B}}_{\smash{k|t}}^{\trs} P\!\right)^{-1} \\
&= (\tilde{E}_{k|t} P + \tilde{F}_{k|t})(\tilde{G}_{k|t} P + \tilde{H}_{k|t})^{-1}
\end{align*}
for $P\in\mathbb{S}_{\smash{++}}^n$, where
\begin{align*}
 &\tilde{E}_{k|t} := \tilde{\lftd{A}}_{\smash{k|t}}^{\trs} + \tilde{\lftd{Q}}_{k|t}\tilde{\lftd{A}}_{\smash{k|t}}^{-1}\tilde{\lftd{B}}_{k|t}\tilde{\lftd{R}}_{\smash{k|t}}^{-1}\tilde{\lftd{B}}_{\smash{k|t}}^{\trs}   \\
 &=\tilde{\lftd{A}}_{\smash{k|t}}^{\trs} ((\tilde{\lftd{B}}_{k|t}\tilde{\lftd{R}}_{\smash{k|t}}^{-1}\tilde{\lftd{B}}_{\smash{k|t}}^{\trs})^{-1} \! +\!(\tilde{\lftd{A}}_{\smash{k|t}}^{\trs})^{-1}\tilde{\lftd{Q}}_{k|t}\tilde{\lftd{A}}_{\smash{k|t}}^{-1}
 ) \tilde{\lftd{B}}_{k|t}\tilde{\lftd{R}}_{\smash{k|t}}^{-1}\tilde{\lftd{B}}_{\smash{k|t}}^{\trs}, 
\end{align*}
$\tilde{F}_{k|t} := \tilde{\lftd{Q}}_{k|t} \tilde{\lftd{A}}_{\smash{k|t}}^{-1}$,
$\tilde{G}_{k|t} := \tilde{\lftd{A}}_{\smash{k|t}}^{-1} \tilde{\lftd{B}}_{k|t} \tilde{\lftd{R}}_{\smash{k|t}}^{-1} \tilde{\lftd{B}}_{\smash{k|t}}^{\trs}$, 
and $\tilde{H}_{k|t} := \tilde{\lftd{A}}_{\smash{k|t}}^{-1}$. Moreover,
\[ ((\tilde{\lftd{B}}_{k|t}\tilde{\lftd{R}}_{\smash{k|t}}^{-1}\tilde{\lftd{B}}_{\smash{k|t}}^{\trs})^{-1} + 
(\tilde{\lftd{A}}_{\smash{k|t}}^{\trs})^{-1}\tilde{\lftd{Q}}_{k|t}\tilde{\lftd{A}}_{\smash{k|t}}^{-1})\in\mathbb{S}_{\smash{++}}^{n}.\]
Therefore, $\tilde{E}_{k|t}$ is non-singular. 
Further,
\begin{align*}
\tilde{F}_{k|t} \tilde{E}_{\smash{k|t}}^{\trs}\! =\! \tilde{\lftd{Q}}_{k|t}\! +\! \tilde{\lftd{Q}}_{k|t}\tilde{\lftd{A}}_{\smash{k|t}}^{-1}\tilde{\lftd{B}}_{k|t}\tilde{\lftd{R}}_{\smash{k|t}}^{-1}\tilde{\lftd{B}}_{\smash{k|t}}^{\trs}(\tilde{\lftd{A}}_{\smash{k|t}}^{\trs})^{-1}\tilde{\lftd{Q}}_{k|t}\!\in
\!\mathbb{S}_{\smash{++}}^n,
\end{align*}
\begin{align*}
&\tilde{E}_{\smash{k|t}}^{\trs}\tilde{G}_{k|t} = 
\tilde{\lftd{B}}_{k|t}\tilde{\lftd{R}}_{\smash{k|t}}^{-1}\tilde{\lftd{B}}_{\smash{k|t}}^{\trs}\nonumber \\
&~+ 
\tilde{\lftd{B}}_{k|t}\tilde{\lftd{R}}_{\smash{k|t}}^{-1}\tilde{\lftd{B}}_{\smash{k|t}}^{\trs}
(\tilde{\lftd{A}}_{\smash{k|t}}^{\trs})^{-1}\tilde{\lftd{Q}}_{k|t}\tilde{\lftd{A}}_{\smash{k|t}}^{-1}
\tilde{\lftd{B}}_{k|t}\tilde{\lftd{R}}_{\smash{k|t}}^{-1}\tilde{\lftd{B}}_{\smash{k|t}}^{\trs} \in\mathbb{S}_{\smash{++}}^n,
\end{align*}
and
\begin{align*} 
&\tilde{G}_{k|t} \tilde{E}_{\smash{k|t}}^{-1} =(\tilde{\lftd{Q}}_{k|t} + \tilde{\lftd{A}}_{k|t}^\trs(\tilde{\lftd{B}}_{k|t} \tilde{\lftd{R}}_{\smash{k|t}}^{-1}\tilde{\lftd{B}}_{k|t}^\trs)^{-1}\tilde{\lftd{A}}_{k|t})^{-1}
\in\mathbb{S}_{\smash{++}}^n,
\end{align*}
since
\begin{align*}
&\tilde{\lftd{A}}_{\smash{k|t}}^{-1} \tilde{\lftd{B}}_{k|t} \tilde{\lftd{R}}_{\smash{k|t}}^{-1}\tilde{\lftd{B}}_{k|t}^\trs
(\tilde{\lftd{A}}_{k|t}^\trs + \tilde{\lftd{Q}}_{k|t} \tilde{\lftd{A}}_{\smash{k|t}}^{-1}\tilde{\lftd{B}}_{k|t}\tilde{\lftd{R}}_{\smash{k|t}}^{-1}\tilde{\lftd{B}}_{k|t}^\trs)^{-1}\\
&\,=\tilde{\lftd{A}}_{\smash{k|t}}^{-1}((\tilde{\lftd{B}}_{k|t} \tilde{\lftd{R}}_{\smash{k|t}}^{-1} \tilde{\lftd{B}}_{k|t}^\trs)^{-1}+(\tilde{\lftd{A}}_{k|t}^\trs)^{-1}\tilde{\lftd{Q}}_{k|t} \tilde{\lftd{A}}_{\smash{k|t}}^{-1})^{-1}(\tilde{\lftd{A}}_{k|t}^\trs)^{-1}.
\end{align*}
As such,~\cite[Thm.~1.7]{bougerol1993kalman} applies to give 
$\delta(\tilde{\lftd{\mathcal{R}}}_{\gamma,k|t}(X), \tilde{\lftd{\mathcal{R}}}_{\gamma,k|t}(P)) \leq \tilde{\rho}_{k|t}\cdot \delta(X,P)$ for all  
$X,P \in \mathbb{S}_{\smash{++}}^n$,
with 
$\tilde{\rho}_{k|t} = \tilde{\zeta}_{k|t}/(\tilde{\zeta}_{k|t} + \tilde{\epsilon}_{k|t})<1$,
$\tilde{\zeta}_{k|t} = 1\big/\lambda_{\min}(
\tilde{F}_{k|t}\tilde{E}_{\smash{k|t}}^{\trs})$,  and
$\tilde{\epsilon}_{k|t} = 
\lambda_{\min}(\tilde{G}_{k|t} \tilde{E}_{\smash{k|t}}^{-1})$, \mc{as per the proof therein}.
\hspace*{\fill}\QEDopen

\subsection{\mc{Riccati operator lifting}} \label{subsec:liftricc}
% Towards lifting~\eqref{eq:RiccatiRecursion} to obtain the Riccati operator~\eqref{eq:lifted_ric_0},
Consider the lifted model formulated in Lemmas~\ref{lemma:lifted_x} and~\ref{lemma:lifted_z}, and corresponding definitions of $\lftd{A}_{k|t}$, $\lftd{B}_{k|t},\lftd{F}_{k|t},\lftd{D}_{k|t},\lftd{E}_{k|t}$, and $\lftd{R}_{k|t}$, \mc{with} $d\in\mathbb{N}$ fixed \mc{according to} Assumption~\ref{asm:uni_obs_ctr}.  Given $t,k\in\mathbb{N}_0$, $\gamma\in\mathbb{R}_{>0}$, 
$\lftd{v}_{k|t} =  (\lftd{u}_{k|t},\lftd{w}_{k|t})
\in\mathbb{R}^{md}\times \mathbb{R}^{nd}$,
$\lftd{x}_{k|t}\in~\mathbb{R}^n$, and $P\in\mathbb{S}^n$, it follows that
%\mc{in view of~\eqref{eq:lifted_dynamic_ch5} and~\eqref{eq:output_z_lifted}},
\begin{align}\label{eq:cost_matrix_quad_lifted}
&\lftd{z}_{\smash{k|t}}^{\trs}\lftd{z}_{k|t} - \gamma^2\lftd{w}_{\smash{k|t}}^{\trs}\lftd{w}_{k|t} + \lftd{x}_{k+1|t}^{\trs} P \lftd{x}_{k+1|t} \nonumber\\
&\quad = \begin{bmatrix} \lftd{x}_{k|t} \\ \lftd{v}_{k|t} \end{bmatrix}^{\trs}
\begin{bmatrix}
\lftd{C}_{\smash{k|t}}^{\trs}\lftd{C}_{k|t}+\lftd{A}_{\smash{k|t}}^{\trs}P\lftd{A}_{k|t} & \lftd{L}_{\smash{k|t}}^{\trs}(P) \\
\lftd{L}_{k|t}(P) & \lftd{M}_{\gamma,k|t}(P)
\end{bmatrix}
\begin{bmatrix} \lftd{x}_{k|t} \\ \lftd{v}_{k|t} \end{bmatrix},
\end{align}
where $\lftd{L}_{k|t}(P) 
:=
\begin{bmatrix}\lftd{D}_{k|t} & \lftd{E}_{k|t}\end{bmatrix}^{\trs}\! \lftd{C}_{k|t} + \tilde{\lftd{B}}_{\smash{k|t}}^{\trs}P\lftd{A}_{k|t}$,
\begin{align}
\label{eq:M_lifted_def}
\lftd{M}_{\gamma,k|t}(P) := 
\tilde{\lftd{R}}_{k|t} + \tilde{\lftd{B}}_{\smash{k|t}}^{\trs} P \tilde{\lftd{B}}_{k|t},
\end{align}
$\tilde{\lftd{B}}_{k|t}$ is defined in~\eqref{eq:btilch5}, and $\gamma$-dependent $\tilde{\lftd{R}}_{k|t}$  \mc{is defined} in~\eqref{eq:rtilch5}.

\begin{lemma} \label{lemma:MhatM}
Given $t\in\mathbb{N}_0$, $\gamma\in\mathbb{R}_{>0}$, and $(\smash{P_s})_{s\in\smash{\{t,\ldots,t+d\}}}\subset\mathbb{S}^n$, suppose for $s\in\{t,\ldots,t+d-1\}$ that $M_{\gamma,s}(P_{s+1})$ as defined in~\eqref{eq:LM_def} is non-singular, and $P_s=\mathcal{R}_{\gamma,s}(P_{s+1})$ \mc{in accordance with}~\eqref{eq:RiccatiGamma}. Then, $\lftd{M}_{\gamma,0|t}(P_{t+d})$ \mc{as defined in}~\eqref{eq:M_lifted_def}
is non-singular.
\end{lemma}
\begin{proof}
If $d=1$, then $\lftd{R}_{0|t} = R_t$, $\lftd{B}_{0|t} = B_t$, $\lftd{F}_{0|t} = I_n$, $\lftd{C}_{0|t} = Q_t^{\smash{1\!/2}}$, $\lftd{D}_{0|t} = 0_{n,m}$, $\lftd{E}_{0|t} = 0_{n,n}$, and $\lftd{F}_{0|t} = I_n$, whereby $\lftd{M}_{\gamma,0|t}(P_{t+1}) = M_{\gamma,t}(P_{t+1})$. 
\mc{As such,} the result follows by induction on noting that it holds with one additional lifting step irrespective of the value of $d$, as shown below. 

Given $\lftd{x}_{0|t}\in\mathbb{R}^n$, $\lftd{v}_{0|t}=\mc{(\lftd{u}_{0|t},\lftd{w}_{0|t})}\in\mathbb{R}^{md}\times\mathbb{R}^{nd}$, $v_{t+d}=(u_{t+d},w_{t+d})\in\mathbb{R}^{m}\times\mathbb{R}^{n}$, and 
$x_{t+d}=\lftd{x}_{1|t} = \lftd{A}_{0|t}\lftd{x}_{0|t} + \tilde{\lftd{B}}_{0|t} \lftd{v}_{0|t}$ as per~\eqref{eq:lifted_dynamic_ch5}, %and~\eqref{eq:btilch5}, 
it follows from~\eqref{eq:ltv_sys_w} that
\begin{align*}
    x_{t+d+1} &= A_{t+d}x_{t+d} + B_{t+d} u_{t+d} + w_{t+d} \\
    &= A_{t+d}\lftd{A}_{0|t}\lftd{x}_{0|t}
    + A_{t+d}\tilde{\lftd{B}}_{0|t} \lftd{v}_{0|t}  
    +
    \!\begin{bmatrix} B_{t+d} & I_n\end{bmatrix} \! v_{t+d}.
\end{align*}
\mc{Furthermore, with $\lftd{z}_{\smash{0|t}}$ as per~\eqref{eq:output_z_lifted}},
% [\begin{smallmatrix}Q_{t+d}^{\smash{1\!/2}} \\ 0_{m,n} \end{smallmatrix}]\lftd{x}_{1|t} + [\begin{smallmatrix}0_{n,m} \\ R_{t+d}^{\smash{1\!/2}}\end{smallmatrix}] u_{t+d}$ 
for any~$P\in\mathbb{S}^n$, 
%\mc{it follows that}
\begin{align}
&\lftd{z}_{\smash{0|t}}^{\trs}\lftd{z}_{0|t} - \gamma^2\lftd{w}_{\smash{0|t}}^{\trs}\lftd{w}_{0|t} \nonumber \\
&+ z_{t+d}^\trs z_{t+d} - \gamma^2w_{t+d}^{\trs}w_{t+d}
+ x_{t+d+1}^{\trs} P x_{t+d+1} \nonumber \\
&= 
\begin{bmatrix} \lftd{x}_{0|t} \\ \lftd{v}_{\smash{0|t}}^{\smash{+}} 
\end{bmatrix}^\trs
\begin{bmatrix}
\lftd{Q}_{\smash{0|t}}^+(P)    & (\lftd{L}_{\smash{0|t}}^{\smash{+}}(P))^{\trs} \\
\lftd{L}_{\smash{0|t}}^{\smash{+}}(P)
& \lftd{M}_{\smash{\gamma,0|t}}^{\smash{+}}(P)
\end{bmatrix}
\begin{bmatrix} \lftd{x}_{0|t} \\ \lftd{v}_{\smash{0|t}}^+ 
\end{bmatrix}, \label{eq:Juwplus}
\end{align}
where ~$z_{t+d} \!=\! 
\begin{bmatrix}
    Q_{t+d}^{\smash{1\!/2}} & 0_{n,m}
\end{bmatrix}^{\trs}\lftd{x}_{1|t} \!+\!
\begin{bmatrix}
   0_{m,n} & R_{t+d}^{\smash{1\!/2}} 
\end{bmatrix}^{\trs}u_{t+d}$ as per~\eqref{eq:output_z_def}, 
$\lftd{v}_{\smash{0|t}}^{\smash{+}} := (\lftd{u}_{0|t}, u_{t+d}, \lftd{w}_{0|t}, w_{t+d})=
\varPi (\lftd{v}_{0|t},v_{t+d})$ for a corresponding permutation matrix $\varPi\in\mathbb{R}^{(m+n)\cdot(d+1) \times (m+n)\cdot(d+1)}$,
\begin{align} 
\lftd{Q}_{\smash{0|t}}^{\smash{+}}(P) &:= \lftd{C}_{\smash{0|t}}^{\trs}\lftd{C}_{0|t} \!+\! \lftd{A}_{\smash{0|t}}^{\trs}(Q_{t+d}\!+\! A_{t+d}^{\trs}PA_{t+d})\lftd{A}_{0|t}, \label{eq:QhatPlus}\\
\label{eq:LhatPlus}
\lftd{L}_{\smash{0|t}}^{\smash{+}}(P) &:= \varPi 
\begin{bmatrix}
\lftd{L}_{\smash{0|t}}^{1+}(P) \\ \lftd{L}_{\smash{0|t}}^{2+}(P)
\end{bmatrix},\\
%\intertext{\mc{where}}
\lftd{L}_{\smash{0|t}}^{1+}(P)&:=\begin{bmatrix}
\lftd{D}_{\smash{0|t}}^{\trs} \\ \lftd{E}_{\smash{0|t}}^{\trs}
\end{bmatrix} \lftd{C}_{0|t}
+ \tilde{\lftd{B}}_{\smash{0|t}}^{\trs} (Q_{t+d}+A_{t+d}^{\trs}PA_{t+d}) \lftd{A}_{0|t} , \label{eq:LhatPlus1}\\
\lftd{L}_{\smash{0|t}}^{2+}(P)
&:=\begin{bmatrix} B_{t+d}^{\trs} \\ I_n \end{bmatrix} P A_{t+d}\lftd{A}_{0|t} = L_{t+d}(P)\lftd{A}_{0|t},
\label{eq:LhatPlus2}
\end{align}
%and $L_{t+d}(P)$ as per~\eqref{eq:LM_def},
and finally,
\begin{align*}
& \lftd{M}_{\mc{\gamma},0|t}^+(P) := \nonumber  \\
&\varPi 
\!\begin{bmatrix}
\tilde{\lftd{R}}_{0|t} \!+\! \tilde{\lftd{B}}_{\smash{0|t}}^\trs (Q_{t+d} \!+\! A_{t+d}^\trs
P A_{t+d}) \tilde{\lftd{B}}_{0|t}
& 
\tilde{\lftd{B}}_{\smash{0|t}}^\trs L_{t+d}^\trs(P) 
\\
L_{t+d}(P)\tilde{\lftd{B}}_{0|t}
& 
M_{\gamma,t+d}(P)
\end{bmatrix}\!
\varPi^{\trs}, 
\end{align*}
with $L_{t+d}(P)$ and $M_{\gamma,t+d}(P)$ as per~\eqref{eq:LM_def}.
\mc{When} $M_{\gamma,t+d}(P)$
is non-singular, Schur decomposition yields
\begin{align}
&\lftd{M}_{\smash{\gamma,0|t}}^+(P) = \nonumber \\
&\varPi \begin{bmatrix} I_{(m+n)d} & (L_{t+d}(P)
\tilde{\lftd{B}}_{0|t}
)^{\trs} (M_{\gamma,t+d}(P))^{-1} \\ 0_{m+n,(m+n)d} & I_{m+n} \end{bmatrix}  
\nonumber \\
& ~\times
\begin{bmatrix}
\lftd{M}_{\gamma,0|t}(\mathcal{R}_{\gamma,t+d}(P))
& 0_{(m+n)d,m+n} \\
0_{m+n,(m+n)d} & M_{\gamma,t+d}(P)
\end{bmatrix}
\nonumber \\
&~
\times \begin{bmatrix} I_{(m+n)d} & 0_{(m+n)d,m+n} \\
(M_{\gamma,t+d}(P))^{-1}L_{t+d}(P) \tilde{\lftd{B}}_{0|t}
& I_{m+n} \end{bmatrix}
\varPi^{\trs}, \label{eq:MhatPlus}
\end{align}
with $\mathcal{R}_{\gamma,t+d}$ and $\lftd{M}_{\gamma,0|t}$ as per~\eqref{eq:RiccatiGamma} and \eqref{eq:M_lifted_def}, respectively.
\mc{Therefore}, if $P=P_{t+d+1}\in\mathbb{S}^n$ with $M_{\gamma,t+d}(P_{t+d+1})$ non-singular, and $P_{t+d} = \mathcal{R}_{\gamma,t+d}(P_{t+d+1})$ with $\lftd{M}_{\gamma,0|t}(P_{t+d})$ non-singular, then $\lftd{M}_{\smash{\gamma,0|t}}^+(P_{t+d+1})$ is non-singular, as required for the induction argument.
\end{proof}

If $\lftd{M}_{\gamma,k|t}(P)$ is non-singular for given $t,k\in\mathbb{N}_0$, $\gamma\in\mathbb{R}_{>0}$, and $P\in\mathbb{S}^n$, then Schur decomposition of~\eqref{eq:cost_matrix_quad_lifted} yields
\begin{align}\label{eq:lifted_ric_fractorization}
&\lftd{z}_{\smash{k|t}}^{\trs}\lftd{z}_{k|t} - \gamma^2\lftd{w}_{\smash{k|t}}^{\trs}\lftd{w}_{k|t} + \lftd{x}_{k+1|t}^{\trs} P \lftd{x}_{k+1|t} 
\nonumber\\
&= \lftd{x}_{\smash{k|t}}^{\trs} \lftd{\mathcal{R}}_{\gamma,k|t}(P)
\lftd{x}_{k|t} \nonumber\\
&\quad + (
\lftd{v}_{k|t}- \lftd{\mathsf{v}}_{\gamma,k|t})^{\trs}\lftd{M}_{\gamma,k|t}(P) (
\lftd{v}_{k|t} - \lftd{\mathsf{v}}_{\gamma,k|t}),
\end{align}
where the lifted $\ell_2$ gain Riccati operator 
\begin{align}
&\lftd{\mathcal{R}}_{\gamma,k|t}(P) := \lftd{C}_{\smash{k|t}}^{\trs}\lftd{C}_{k|t}+\lftd{A}_{\smash{k|t}}^{\trs}P\lftd{A}_{k|t}\nonumber \\
& \qquad \qquad \qquad - \lftd{L}_{\smash{k|t}}^{\trs}(P) (\lftd{M}_{\gamma,k|t}(P))^{-1} \lftd{L}_{k|t}(P), \label{eq:lifted_ric}
\end{align}
and $\lftd{\mathsf{v}}_{\gamma,k|t}:= -(\lftd{M}_{\gamma,k|t}(P))^{-1}\lftd{L}_{k|t}(P)\lftd{x}_{k|t}$.

\begin{lemma}\label{lemma:ric_gamma_comp}
Given  
$t\in\mathbb{N}_0$, $\gamma\in\mathbb{R}_{>0}$, and $(\smash{P_s})_{s\in\smash{\{t,\ldots,t+d\}}}\subset\mathbb{S}^n$, suppose for $s\in\{t,\ldots,t+d-1\}$ that $M_{\gamma,s}(P_{s+1})$ 
as defined in~\eqref{eq:LM_def} is non-singular, and $P_s=\mathcal{R}_{\gamma,s}(P_{s+1})$ \mc{in accordance with}~\eqref{eq:RiccatiGamma}.
Then,
$\lftd{\mathcal{R}}_{\gamma,0|t}(P_{t+d}) = \mathcal{R}_{\gamma,t} \circ \cdots \circ \mathcal{R}_{\gamma,t+d-1}(P_{t+d})=P_t$.
\end{lemma}
\begin{proof}
Building upon the induction argument in the proof of Lemma~\ref{lemma:MhatM} yields the result. If $d=1$, then $\lftd{\mathcal{R}}_{\gamma,0|t}(P_{t+1}) = \mathcal{R}_{\gamma,t}(P_{t+1})$ by definition. Given this, it remains to show that if $\lftd{\mathcal{R}}_{\gamma,0|t}(P_{t+d})
=\mathcal{R}_{\gamma,t} \circ \cdots \circ \mathcal{R}_{\gamma,t+d-1}(P_{t+d})$, then it is also true with one additional step in the lifting. To this end, observe that for $P\in\mathbb{S}^n$ such that $\lftd{M}_{\smash{\gamma,0|t}}^+(P)$ in~\eqref{eq:MhatPlus} is non-singular, the lifted Riccati operator with one additional step is given by the corresponding Schur complement
\begin{align} \label{eq:RiccPlus}
\lftd{\mathcal{R}}_{\smash{\gamma,0|t}}^{\smash{+}}(P) & :=
\lftd{Q}_{\smash{0|t}}^+(P) 
- (\lftd{L}_{\smash{0|t}}^{\smash{+}}(P))^\trs     
(\lftd{M}_{\smash{\gamma,0|t}}^+(P))^{-1}
\lftd{L}_{\smash{0|t}}^{\smash{+}}(P)
\end{align}
of~\eqref{eq:Juwplus}. Noting that
\begin{align*}
&(\lftd{M}_{\smash{\gamma,0|t}}^+(P))^{-1} =   \\
&
\varPi 
\begin{bmatrix} I_{(m+n)d} & 0_{(m+n)d,m+n} 
\\
-(M_{\gamma,t+d}(P))^{-1}L_{t+d}(P) \tilde{\lftd{B}}_{0|t}
& I_{m+n} \end{bmatrix} 
\\
& ~ \times
\begin{bmatrix}
(\lftd{M}_{\gamma,0|t}(\mathcal{R}_{\gamma,t+d}(P)))^{-1}
& 0_{(m+n)d,m+n} 
\\
0_{m+n,(m+n)d} & (M_{\gamma,t+d}(P))^{-1}
\end{bmatrix}
 \\
& 
~ \times
\begin{bmatrix} I_{(m+n)d} & -((M_{\gamma,t+d}(P))^{-1}L_{t+d}(P) \tilde{\lftd{B}}_{0|t})^{\trs} \\ 0_{m+n,(m+n)d}
& I_{m+n} \end{bmatrix}
\varPi^{\trs},
\end{align*}
%if
%$M_{\gamma,t+d}(P)$ and $\lftd{M}_{\gamma,0|t}(P)$ as per~\eqref{eq:LM_def} and~\eqref{eq:M_lifted_def} are non-singular,
it follow from~\eqref{eq:RiccatiGamma},~\eqref{eq:QhatPlus}--\eqref{eq:LhatPlus2},~\eqref{eq:lifted_ric}, and~\eqref{eq:RiccPlus}, that
\begin{align*}
&\lftd{\mathcal{R}}_{\smash{\gamma,0|t}}^{\smash{+}}(P)\\
&=
\lftd{C}_{\smash{0|t}}^\trs \lftd{C}_{0|t} +
\lftd{A}_{\smash{0|t}}^\trs \mathcal{R}_{\gamma,t+d}(P)
\lftd{A}_{0|t} \\
&
%\quad 
- \lftd{L}_{0|t}(\mathcal{R}_{\gamma,t+d}(P))^\trs (\lftd{M}_{\gamma,0|t}(\mathcal{R}_{\gamma,t+d}(P)))^{-1}
\lftd{L}_{0|t}(\mathcal{R}_{\gamma,t+d}(P))\\
& =
\lftd{\mathcal{R}}_{\gamma,0|t}(\mathcal{R}_{\gamma,t+d}(P)).
\end{align*}
So, if $P=P_{t+d+1}\in\mathbb{S}^n$ with $M_{\gamma,t+d}(P_{t+d+1})$ non-singular, and $P_{t+d}=\mathcal{R}_{\gamma,t+d}(P_{t+d+1})$, then $\lftd{\mathcal{R}}_{\smash{\gamma,0|t}}^+(P_{t+d+1}) = \lftd{\mathcal{R}}_{\gamma,0|t}(P_{t+d})=\mathcal{R}_{\gamma,t} \circ \cdots \circ \mathcal{R}_{\gamma,t+d-1}(P_{t+d})=\mathcal{R}_{\gamma,t} \circ \cdots \circ \mathcal{R}_{\gamma,t+d}(P_{t+d+1})$, as required for the induction argument.
\end{proof}

\begin{remark}\label{rem:lift_equiv}
    As noted in Remark~\ref{rem:baseline}, with $(P_t)_{t\in\mathbb{N}_0}\subset\mathbb{S}_{\smash{+}}^n$ \mc{satisfying} part 1) of the hypothesis in Theorem~\ref{theorem:h_inf_T_lb}, $M_{\gamma,t}(P_{t+1})$ is non-singular, and $P_t=\mathcal{R}_{\gamma,t}(P_{t+1})$ for every $t\in\mathbb{N}_0$. As such, \mc{for all $t,k\in\mathbb{N}_0$,} $\lftd{M}_{\gamma,0|t}(P_{t+d})$ and $\lftd{M}_{\gamma,k|t}(P_{t+d(k+1)})$ are non-singular by Lemma~\ref{lemma:MhatM}, since $\lftd{M}_{\gamma,k|t}(P_{t+d(k+1)})=\lftd{M}_{\gamma,0|t+dk}(P_{t+dk+d})$. Further,  
    $P_{t+dk}=\mathcal{R}_{\gamma,t+dk} \circ \cdots \circ \mathcal{R}_{\gamma,t+d(k+1)-1}(P_{t+d(k+1)})=\lftd{\mathcal{R}}_{\gamma,k|t}(P_{t+d(k+1)})$ by Lemma~\ref{lemma:ric_gamma_comp}, \mc{as exploited subsequently.}
\end{remark}
\begin{remark} \label{rem:prePDP}
Given $t,k\in\mathbb{N}_0$, suppose $\gamma\in\mathbb{R}_{>0}$ and $P\in\mathbb{S}_{\smash{+}}^{n}$ are such that $\lftd{M}_{\gamma,k|t}(P)$ in~\eqref{eq:M_lifted_def} 
is non-singular. Then,
%from~\eqref{eq:lifted_ric},
\begin{align*}
&\lftd{\mathcal{R}}_{\gamma,k|t}(P) \nonumber \\
&= \lftd{C}_{\smash{k|t}}^{\trs} \lftd{C}_{k|t} + \lftd{A}_{\smash{k|t}}^{\trs} P \lftd{A}_{k|t}  \nonumber \\
& \quad - 
\left(\! \lftd{A}_{\smash{k|t}}^{\trs}P
\tilde{\lftd{B}}_{k|t} 
+ \lftd{C}_{\smash{k|t}}^{\trs} \!\begin{bmatrix} \lftd{D}_{k|t} & \lftd{E}_{k|t} \end{bmatrix} \!\right) \nonumber \\
&\quad\quad\times
(\tilde{\lftd{R}}_{k|t} + \tilde{\lftd{B}}_{k|t}^\trs P \tilde{\lftd{B}}_{k|t})^{-1}
\left( \!
\tilde{\lftd{B}}_{k|t}^\trs 
 P \lftd{A}_{k|t} \!+\!\! \begin{bmatrix} \lftd{D}_{\smash{k|t}}^{\trs} \\ \lftd{E}_{\smash{k|t}}^{\trs} \end{bmatrix} \! \lftd{C}_{k|t} \!\right).
\end{align*}  
With $\tilde{\lftd{R}}_{k|t}\in\mathbb{S}^{(m+n)d}$ non-singular, it follows that 
%with reference to
$\tilde{\lftd{\mathcal{R}}}_{\gamma,k|t}(P)=\lftd{\mathcal{R}}_{\gamma,k|t}(P)$  in~\eqref{eq:lifted_ric_0} by~\cite[Prop.~12.1.1]{lancaster1995algebraic}. \mc{Further,} 
\begin{align} \label{eq:transformedRicc2}
&\tilde{\lftd{\mathcal{R}}}_{\gamma,k|t}(P)
= \nonumber \\
&\tilde{\lftd{Q}}_{k|t}
+ \tilde{\lftd{A}}_{k|t}^\trs P^{\smash{1\!/2}}
(I + P^{\smash{1\!/2}}\tilde{\lftd{B}}_{k|t}\tilde{\lftd{R}}_{\smash{k|t}}^{-1} \tilde{\lftd{B}}_{\smash{k|t}}^{\trs}P^{\smash{1\!/2}})^{-1}
P^{\smash{1\!/2}} \tilde{\lftd{A}}_{k|t}
\end{align}
by the Woodbury formula, \mc{whenever} $\tilde{\lftd{B}}_{k|t}\tilde{\lftd{R}}_{\smash{k|t}}^{-1} \tilde{\lftd{B}}_{\smash{k|t}}^{\trs}\in\mathbb{S}_{\smash{\smash{++}}}^n$.
As such, \mc{if} $\tilde{\lftd{Q}}_{k|t}\in\mathbb{S}_{\smash{++}}^{(m+n)d}$ also holds, \mc{then} $\tilde{\lftd{\mathcal{R}}}_{\gamma,k|t}(P)\in\mathbb{S}_{\smash{++}}^n$. Finally, 
from~\eqref{eq:transformedRicc2},
\begin{align}\label{eq:lifted_ric_2}
\tilde{\lftd{\mathcal{R}}}_{\gamma,k|t}(P) 
&= \tilde{\lftd{Q}}_{k|t} + \tilde{\lftd{A}}_{\smash{k|t}}^{\trs}(P^{-1} + \tilde{\lftd{B}}_{k|t}\tilde{\lftd{R}}_{\smash{k|t}}^{-1}\tilde{\lftd{B}}_{\smash{k|t}}^{\trs} )^{-1}\tilde{\lftd{A}}_{k|t}
\end{align}
for $P\in\mathbb{S}_{\smash{++}}^n$.
\end{remark}
\begin{remark} \label{rem:PDP}
As noted in Remark~\ref{rem:lift_equiv}, with 
$(P_t)_{t\in\mathbb{N}_0}\subset\mathbb{S}_{\smash{+}}^n$ as per part 1) of the hypothesis in Theorem~\ref{theorem:h_inf_T_lb}, $\lftd{M}_{\gamma,k|t}(P_{t+d(k+1)})$ is non-singular for all $t,k\in\mathbb{N}_0$. With part 2) of the hypothesis, 
$\tilde{\lftd{R}}_{k|t}=\tilde{\lftd{R}}_{0|t+kd}$ is non-singular,
$\tilde{\lftd{Q}}_{k|t}=\tilde{\lftd{Q}}_{0|t+kd}\in\mathbb{S}_{\smash{++}}^{n}$, and $\tilde{\lftd{B}}_{k|t}\tilde{\lftd{R}}_{\smash{k|t}}^{-1} \tilde{\lftd{B}}_{\smash{k|t}}^{\trs} = \tilde{\lftd{B}}_{0|t+kd}\tilde{\lftd{R}}_{\smash{0|t+kd}}^{-1} \tilde{\lftd{B}}_{0|t+kd}^{\trs}\in\mathbb{S}_{\smash{++}}^{n}$. Thus, in view of Remark~\ref{rem:prePDP}, $P_{t+dk}=\lftd{\mathcal{R}}_{\gamma,k|t}(P_{t+d(k+1)})\in\mathbb{S}_{\smash{++}}^n$; \mc{i.e.,} 
the hypothesized 
$(P_t)_{t\in\mathbb{N}_0}\subset\mathbb{S}_{\smash{+}}^n$ %in Theorem~\ref{theorem:h_inf_T_lb} is 
\mc{must be} contained in $\mathbb{S}_{\smash{++}}^n$ \mc{as previously noted in Remark~\ref{rem:baseline}}.
\end{remark}

\begin{remark} \label{rem:part2forallt}
    With part 2) of the hypothesis in Theorem~\ref{theorem:h_inf_T_lb},  $\tilde{\lftd{R}}_{k|t}=\tilde{\lftd{R}}_{0|t+dk}$ is non-singular for all $t,k\in\mathbb{N}_0$, and thus, $\tilde{\lftd{A}}_{k|t}=\tilde{\lftd{A}}_{0|t+dk}$ \mc{is} non-singular by Lemma~\ref{lemma:lifted_A_nonsingular} \mc{in the Appendix}. Furthermore, $\tilde{\lftd{B}}_{k|t}\tilde{\lftd{R}}_{\smash{k|t}}^{-1}\tilde{\lftd{B}}_{\smash{k|t}}^{\trs}=\tilde{\lftd{B}}_{0|t+dk}\tilde{\lftd{R}}_{0|t+dk}^{-1}\tilde{\lftd{B}}_{0|t+dk}^{\trs} \in \mathbb{S}_{\smash{++}}^n$, and $\tilde{\lftd{Q}}_{k|t}=\tilde{\lftd{Q}}_{0|t+dk} \in \mathbb{S}_{\smash{++}}^n$. Therefore,
    $\tilde{\lftd{\mathcal{R}}}_{\gamma,k|t}(X)\in\mathbb{S}_{\smash{++}}^n$
over $X\in\mathbb{S}_{\smash{++}}^n$ by Remark~\ref{rem:prePDP}, and 
$\tilde{\lftd{\mathcal{R}}}_{\gamma,k|t}(\cdot)$ is a strict contraction by~Theorem~\ref{theorem:ric_gamma_contraction}. 
\end{remark}

\subsection{\mc{Continuity of closed-loop} \texorpdfstring{$\ell_2$}{l2} gain \mc{performance}} 
\label{subsec:approx}

% In Theorem~\ref{theorem:delta_requirement}, the hypothesis~\eqref{eq:X_requirement_final} implies the control policy~\eqref{eq:h_inf_approx_policy} meets the given performance loss specification $\beta\in\mathbb{R}_{>0}$, relative to the baseline gain bound $\gamma\in\mathbb{R}_{>0}$. In Section~\ref{subsec:proof}, where Theorem~\ref{theorem:h_inf_T_lb} is proved, it is established that the construction of  $(X_{t+1})_{t\in\mathbb{N}_0}$ in~\eqref{eq:h_inf_approx_X} satisfies this hypothesis. 
The following lemmas lead to a proof of Theorem~\ref{theorem:delta_requirement}.

\begin{lemma}\label{lemma:l2_sufficient_1}
Given $\gamma\in\mathbb{R}_{>0}$, suppose $(P_t)_{t\in\mathbb{N}_0}\subset\mathbb{S}_{\smash{++}}^n$ is such that~\eqref{eq:sc_P_ub} and~\eqref{eq:RiccatiRecursion} hold. Further, given $\beta\in\mathbb{R}_{>0}$, and bounded sequence $(X_{t+1})_{t\in\mathbb{N}_0}\subset\mathbb{S}_{\smash{++}}^n$, suppose 
\begin{align*}
&(\exists \varepsilon\in\mathbb{R}_{>0})~(\forall k \in \mathbb{N}_0)\\
&\quad (\,X_{t+1} -(\gamma+\beta)^2 I_n \preceq -\varepsilon I_n\,) \\
&\quad~ \land ~ 
  (\,P_{t+1}\preceq X_{t+1}\,) ~ \land ~ (\,\mathcal{R}_{\gamma+\beta,t}(X_{t+1}) \preceq \mathcal{R}_{\gamma,t}(P_{t+1})\,),
\end{align*}
in accordance with~\eqref{eq:RiccatiGamma}.
Then, with $J_{\gamma+\beta}$ as per~\eqref{eq:J_def}, the state feedback control policy $u_t\!=\!\mathsf{u}_{\smash{\gamma+\beta,t}}^{\smash{\mathsf{fin}}}(x_t)$ defined in~\eqref{eq:h_inf_approx_policy} for the system~\eqref{eq:ltv_sys_w}--\eqref{eq:output_z_def} achieves 
$(\forall w\in\ell_2)~J_{\gamma+\beta}(u,w) \leq 0$.
\end{lemma}

\begin{proof}
Since $((\gamma+\beta)^2 I_n - X_{t+1}) \in \mathbb{S}_{\smash{++}}^n$, and $R_t\in\mathbb{S}_{\smash{++}}^n$,
\begin{align*}
&\nabla_{\!\gamma+\beta,t}(X_{t+1}) \\
&= R_t  +
B_t^{\trs}(X_{t+1} + X_{t+1}((\gamma+\beta)^2I_n - X_{t+1})^{-1}X_{t+1})B_t \\
&\in\mathbb{S}_{\smash{++}}^m,
\end{align*}
This is the Schur complement of 
$M_{\gamma+\beta,t}(X_{t+1})$, as defined in~\eqref{eq:LM_def}, which is therefore non-singular. Given $x_t \in \mathbb{R}^n$, and $v_{t}=(u_t,w_t) \in \mathbb{R}^{m+n}$, it follows from~\eqref{eq:ltv_sys_w}--\eqref{eq:output_z_def} that
\begin{align*}
&z_t^{\trs}z_t - (\gamma+\beta)^2w_t^{\trs}w_t + x_{t+1}^{\trs}X_{t+1}x_{t+1} \nonumber\\
&=x_{t}^{\trs} \mathcal{R}_{\gamma+\beta,t}(X_{t+1}) x_{t} \nonumber \\
&\quad + (
v_{t}\!-\! \mathsf{v}_{\smash{\gamma+\beta,t}}^{\smash{\mathsf{fin}}})^{\trs} M_{\gamma+\beta,t}(X_{t+1}) (
 v_{t}\! -\! \mathsf{v}_{\smash{\gamma+\beta,t}}^{\smash{\mathsf{fin}}}), 
\end{align*}
with 
$\mathsf{v}_{\smash{\gamma+\beta,t}}^{\smash{\mathsf{fin}}} := -(M_{\gamma+\beta,t}(X_{t+1}))^{-1}L_t(X_{t+1}) x_t$,
and $L_t(X_{t+1})$ as per~\eqref{eq:LM_def}; cf.~derivation of~\eqref{eq:lifted_ric_fractorization} from~\eqref{eq:cost_matrix_quad_lifted} by Schur decomposition. Then, further
decomposition of $M_{\gamma+\beta,t}(X_{t+1})$ in the same fashion leads to
\begin{align}\label{eq:cost_complete_square_ch5}
&z_t^{\trs}z_t - (\gamma+\beta)^2w_t^{\trs}w_t + x_{t+1}^{\trs}X_{t+1}x_{t+1} \nonumber\\
&= x_t^{\trs} \mathcal{R}_{\gamma+\beta,t}(X_{t+1}) x_t \nonumber \\
& \quad + (u_t-\mathsf{u}_{\smash{\gamma+\beta,t}}^{\smash{\mathsf{fin}}})^{\trs} \nabla_{\!\gamma+\beta,t}(X_{t+1}) (u_t-\mathsf{u}_{\smash{\gamma+\beta,t}}^{\smash{\mathsf{fin}}}) \nonumber\\
& \quad + (w_t-\mathsf{w}_{\smash{\gamma+\beta,t}}^{\smash{\mathsf{fin}}})^{\trs}(X_{t+1}-(\gamma+\beta)^2I_n)(w_t\!-\!\mathsf{w}_{\smash{\gamma+\beta,t}}^{\smash{\mathsf{fin}}}),
\end{align}
with 
$\mathsf{w}_{\smash{\gamma+\beta,t}}^{\smash{\mathsf{fin}}} = -(X_{t+1} - (\gamma+\beta)^2I_n)^{-1}X_{t+1}(A_tx_t+B_tu_t)$, and $\mathsf{u}_{\smash{\gamma+\beta,t}}^{\smash{\mathsf{fin}}}$ as per~\eqref{eq:h_inf_approx_policy}.

Now given any $w = (w_t)_{t\in\mathbb{N}_0}\in\ell_2$, with $u_t=\mathsf{u}_{\smash{\gamma+\beta,t}}^{\smash{\mathsf{fin}}}(x_t)$ in~\eqref{eq:ltv_sys_w}--\eqref{eq:output_z_def}, 
\begin{align*}
&\sum_{t=0}^N z_t^{\trs}z_t - (\gamma+\beta)^2 w_t^{\trs}w_t \nonumber\\
&\leq \sum_{t=0}^N x_t^{\trs} \mathcal{R}_{\gamma+\beta,t}(X_{t+1}) x_t - x_{t+1}^{\trs}X_{t+1}x_{t+1} \nonumber\\
&\leq \sum_{t=0}^N \left( x_t^{\trs} \mathcal{R}_{\gamma+\beta,t}(X_{t+1}) x_t - x_t^{\trs}\mathcal{R}_{\gamma,t}(P_{t+1})x_t\right. \nonumber \\
& \qquad \qquad \qquad \qquad \qquad\quad  \left. + x_t^{\trs}P_tx_t - x_{t+1}^{\trs}P_{t+1}x_{t+1}\right) \nonumber\\
&= x_0^{\trs}P_0x_0 - x_{N+1}^{\trs}P_{N+1}x_{N+1} \nonumber \\
&\qquad \qquad + \sum_{t=0}^N x_t^{\trs} \left( \mathcal{R}_{\gamma+\beta,t}(X_{t+1}) - \mathcal{R}_{\gamma,t}(P_{t+1}) \right) x_t \nonumber\\
&\leq \sum_{t=0}^N x_t^{\trs} \left( \mathcal{R}_{\gamma+\beta,t}(X_{t+1}) - \mathcal{R}_{\gamma,t}(P_{t+1}) \right) x_t,
\end{align*}
for every $N\in\mathbb{N}$.
The first inequality follows from \eqref{eq:cost_complete_square_ch5} since $X_{t+1} -(\gamma+\beta)^2 I_n \preceq -\varepsilon I_n$.
%$(\gamma+\beta)^2 I_n - X_{t+1} \in \mathbb{S}_{\smash{++}}^n$. 
The second inequality follows from $P_t = \mathcal{R}_{\gamma,t}(P_{t+1})$ as per~\eqref{eq:RiccatiRecursion}, and the hypothesis $P_{t+1}\preceq X_{t+1}$. The last inequality holds because $x_0 = 0$, and $P_{N+1} \in \mathbb{S}_{\smash{++}}^n$. 
Therefore, with the hypothesis $\mathcal{R}_{\gamma+\beta,t}(X_{t+1}) \preceq \mathcal{R}_{\gamma,t}(P_{t+1})$, 
\begin{align*}
    (\forall N\in\mathbb{N}) \quad
    \sum_{t=0}^{N-1} z_t^{\trs}z_t - (\gamma+\beta)^2 w_t^{\trs}w_t \leq 0.
\end{align*}
As such, $z \in \ell_2$ since $w \in \ell_2$, with $\|z\|_2^2 \leq (\gamma+\beta)^2\|w\|_2^2$.
\end{proof}

\begin{lemma}\label{lemma:ric_gamma_new}
% Under Assumption~\ref{asm:sc}, 
Given $\gamma\!\in\!\mathbb{R}_{>0}$, and $P\! \in\! \mathbb{S}_{\smash{++}}^n$ such that $P \!\prec\! \gamma^2 I_n$, the $\ell_2$ gain Riccati operator can be expressed as
\begin{align} \label{eq:RiccAlt}
\mathcal{R}_{\gamma,t}(P) = Q_t + A_t^{\trs}(P^{-1} - \gamma^{-2}I_n + B_tR_t^{-1}B_t^{\trs} )^{-1}A_t
\end{align}
in accordance with~\eqref{eq:RiccatiGamma}.
\end{lemma}

\begin{proof}
Since $(P-\gamma^2 I_n)$ is non-singular by hypothesis, Schur decomposition of $M_{\gamma,t}(P)$, as defined in~\eqref{eq:LM_def}, yields
\begin{align}\label{eq:proof_rewrite_ric}
&\mathcal{R}_{\gamma,t}(P) \nonumber \\
&= Q_t + A_t^{\trs}PA_t \nonumber \\
&~~~ -\! 
\begin{bmatrix}
B_t^{\trs}(P - P(P-\gamma^2I_n)^{-1}P)A_t \\
PA_t
\end{bmatrix}^{\trs} \nonumber \\
&\qquad\qquad\times
\begin{bmatrix}
(\nabla_{\!\gamma,t}(P))^{-1} & 0_{m,n} \\
0_{n,m} & \! (P\!-\!\gamma^2I_n)^{-1} \!
\end{bmatrix}\nonumber \\
&\qquad\qquad\qquad \times
\begin{bmatrix}
B_t^{\trs}(P - P(P-\gamma^2I_n)^{-1}P)A_t \\
PA_t
\end{bmatrix},
\end{align}
where $\nabla_{\!\gamma,t}(P) = R_t + B_t^{\trs}(P - P(P-\gamma^2 I_n)^{-1}P)B_t$. Thus, with $W:=P-P(P-\gamma^2 I_n)^{-1}P$, 
\begin{align*}
\mathcal{R}_{\gamma,t}(P)
&= Q_t + A_t^{\trs}(W 
-WB_t  
(R_t + B_t^{\trs}WB_t)^{-1}B_t^{\trs}W)A_t.
\end{align*}
As such, the expression~\eqref{eq:RiccAlt} for $\mathcal{R}_{\gamma,t}(P)$ follows from the identities $W=~\!\!(P^{-1}\!-\!\gamma^{-2}I_n)^{-1}$ and $(W-WB_t(R_t + B_t^{\trs}WB_t)^{-1}B_t^{\trs}W)=(W^{-1}+B_tR_t^{-1}B_t)^{-1}$, which both hold by the Woodbury formula. 
\end{proof}

\begin{lemma}\label{lemma:X_requirement} 
Given $\gamma, \beta \in \mathbb{R}_{>0}$, and $P, X \in \mathbb{S}_{\smash{++}}^n$, if 
\begin{align}\label{eq:X_requirement}
(P \prec \gamma^2 I_n) ~ \mc{\land} ~
(P^{-1} \!-\! X^{-1} \prec \left( \gamma^{-2} \!-\! (\gamma+\beta)^{-2} \right) I_n),
\end{align}
then $X \prec (\gamma+\beta)^2 I_n$ and $\mc{(\forall t\in\mathbb{N}_0)}~\mathcal{R}_{\gamma+\beta,t}(X) \preceq \mathcal{R}_{\gamma,t}(P)$.
%for all $t \in \mathbb{N}_0$.
\end{lemma}

\begin{proof}
From~\eqref{eq:X_requirement}, 
\begin{align}\label{eq:x_p_inv_compare}
0 \prec P^{-1} - \gamma^{-2}I_n \prec X^{-1} - (\gamma+\beta)^{-2}I_n,
\end{align}
which implies $X \prec (\gamma+\beta)^2 I_n$.
% because $(Y^{-1}-Z^{-1})^{-1} = Y + Y(Z-Y)^{-1}Y$ and $(Z-Y)^{-1} = Z^{-1} + Z^{-1}(Y^{-1}-Z^{-1})^{-1}Z^{-1}$ for all $Y,Z\in\mathbb{S}_{\smash{++}}^n$ by the Woodbury formula, and therfore, $0\prec Y\prec Z \Leftrightarrow 0\prec Z^{-1} \prec Y^{-1}$. 
Given $B_t R_t^{-1} B_t^\trs \in \mathbb{S}_{\smash{+}}^n$, it also follows from~\eqref{eq:x_p_inv_compare} that
\begin{align*}
&\left( X^{-1} - (\gamma+\beta)^{-2}I_n + B_tR_t^{-1}B_t^{\trs} \right)^{-1}\\
&\qquad\qquad\qquad \prec \left( P^{-1} \!-\! \gamma^{-2} I_n + B_tR_t^{-1}B_t^{\trs} \right)^{-1}.
\end{align*} 
As such, 
\begin{align*}
&Q_t + A_t^{\trs} \left( X^{-1} \!-\! (\gamma+\beta)^{-2}I_n \!+\! B_tR_t^{-1}B_t^{\trs} \right)^{-1} A_t \\
&\qquad\qquad\qquad \preceq 
Q_t + A_t^{\trs} \left( P^{-1} \!-\! \gamma^{-2} I_n \!+\! B_tR_t^{-1}B_t^{\trs} \right)^{-1} A_t,
\end{align*}
which is $\mathcal{R}_{\gamma+\beta,t}(X) \preceq \mathcal{R}_{\gamma,t}(P)$ in view of Lemma~\ref{lemma:ric_gamma_new}.
\end{proof}

{\noindent\hspace{1em}{\bf\itshape Proof of  Theorem~\ref{theorem:delta_requirement}:}
% %Since 
% First note that $0\prec Y\prec Z \Leftrightarrow 0\prec Z^{-1} \prec Y^{-1}$,
% %for all $Y,Z\in\mathbb{S}_{\smash{++}}^n$, 
% given that $(Y^{-1}-Z^{-1})^{-1} = Y + Y(Z-Y)^{-1}Y$ and $(Z-Y)^{-1} = Z^{-1} + Z^{-1}(Y^{-1}-Z^{-1})^{-1}Z^{-1}$ 
% %for all $Y,Z\in\mathbb{S}_{\smash{++}}^n$, 
% by the Woodbury formula.
%as noted below~\eqref{eq:x_p_inv_compare}, and since
Since $0\prec P_{t+1} \prec X_{t+1}$ by hypothesis, it follows that $0\prec X_{t+1}^{-1} \prec P_{t+1}^{-1}$, \mc{which yields}
\begin{align} \label{eq:lambdamin}
    \lambda_{\max}(X_{t+1}^{-1})  < \lambda_{\max}(P_{t+1}^{-1}) = 1/\lambda_{\min}(P_{t+1}) \leq 1\big/\underline{\lambda}
\end{align}
and $I_n \prec X_{t+1}^{\smash{1\!/2}} P_{t+1}^{-1} X_{t+1}^{\smash{1\!/2}}$, \mc{whereby} $\{\lambda_1,\ldots,\lambda_n\}=\mathrm{spec}(X_{t+1}^{\smash{1\!/2}}P_{t+1}^{-1}X_{t+1}^{\smash{1\!/2}}) 
=\mathrm{spec}(X_{t+1}P_{t+1}^{-1})
\subset\mathbb{R}_{>1}$. As such, $\log(\lambda_i)>0$ for all $i\in\{1,\ldots,n\}$, \mc{and} 
\begin{align}
\delta(X_{t+1},P_{t+1})
&={\textstyle \sqrt{\sum_{i\in\{1,\ldots,n\}} \big(\log(\lambda_i)\big)^2} } \nonumber \\
&\geq {\textstyle \sqrt{\max_{i\in\{1,\ldots,n\}} \big(\log(\lambda_i)\big)^2} } \nonumber \\
&= {\textstyle \sqrt{\big(\max_{i\in\{1,\ldots,n\}} \log( \lambda_i)\big)^2}} \nonumber \\
&= \log(\lambda_{\max}(X_{t+1}^{\smash{1\!/2}}P_{\smash{t+1}}^{-1}X_{t+1}^{\smash{1\!/2}}))>0
\label{eq:deltalambdamax}
\end{align}
\mc{because} $\lambda \mapsto \lambda^2$ and $\lambda \mapsto \log(\lambda)$ are both increasing functions over $\lambda\in\mathbb{R}_{>0}$.

\mc{Given} $\lambda_{\max}(YZ) \leq \lambda_{\max}(Y) \lambda_{\max}(Z)$ and $\lambda_{\max}(Z^{\smash{1\!/2}})=\sqrt{\lambda_{\max}(Z)}$ for all $Y,Z\in\mathbb{S}_{\smash{++}}^n$, it follows from~\eqref{eq:deltalambdamax},~\eqref{eq:lambdamin},~\eqref{eq:alpha_ch5}, and~\eqref{eq:X_requirement_final}, that 
\begin{align*}
&\lambda_{\max}(P_{t+1}^{-1} - X_{t+1}^{-1}) \\
&\leq \lambda_{\max}(X_{t+1}^{-\smash{1\!/2}}) \lambda_{\max}(X_{t+1}^{\smash{1\!/2}}P_{t+1}^{-1}X_{t+1}^{\smash{1\!/2}} - I_n)
\lambda_{\max}(X_{t+1}^{-\smash{1\!/2}})\\
& \leq 
\lambda_{\max}(X_{t+1}^{-1}) \big( \exp(\delta(X_{t+1}, P_{t+1})) - 1 \big) \\
&\leq \frac{1}{\underline{\lambda}} 
\big( \exp(\delta(X_{t+1}, P_{t+1})) - 1 \big) \\
&\leq (\gamma^{-2} - (\gamma+\beta)^{-2}) - \varepsilon.
\end{align*}
 Therefore, \mc{uniformity of the} hypothesis $P_{t+1}- \gamma^2 I_n \prec -\varepsilon I_n$ in~\eqref{eq:sc_P_ub} \mc{means that} both parts of the sufficient \mc{conjunction}~\eqref{eq:X_requirement} in Lemma~\ref{lemma:X_requirement} hold with some margin. \mc{Thus,} $\mathcal{R}_{\gamma+\beta,t}(X_{t+1})\preceq \mathcal{R}_{\gamma,t}(P_{t+1})$ and $X_{t+1}-(\gamma+\beta)^2 I_n \prec 0$ uniformly. Lemma~\ref{lemma:l2_sufficient_1} then applies \mc{to yield} the stated closed-loop $\ell_2$ gain bound. 
\hspace*{\fill}\QEDopen

\subsection{Proof of Theorem~\ref{theorem:h_inf_T_lb}} \label{subsec:proof}

To complete the proof, the context \mc{here} returns to the lifted domain of Sections~\ref{subsec:strcontract} and~\ref{subsec:liftricc}, starting with a collection of lemmas, which \mc{for convenience} all \mc{refer to} the following summary of the \mc{two parts of the hypothesis on the baseline gain bound $\gamma\in\mathbb{R}_{>0}$} in Theorem~\ref{theorem:h_inf_T_lb}. 
\begin{hypothesis} \label{asm:gamprf}
Part 1) and \mc{Part} 2) of the hypothesis in Theorem~\ref{theorem:h_inf_T_lb} hold, with $\tilde{B}_{k|t}=\tilde{B}_{0|t+kd}$, $\tilde{R}_{k|t}=\tilde{R}_{0|t+kd}$,  $\tilde{Q}_{k|t}=\tilde{Q}_{0|t+kd}$, and $\tilde{A}_{k|t}=\tilde{A}_{0|t+kd}$
% where
% %by~\eqref{eq:sc_P_ub} 
% %and~\eqref{eq:RiccatiRecursion} for some 
% %$(P_t)_{t\in\mathbb{N}_0}\subset\mathbb{S}_{\smash{++}}^n$ by Remark~\ref{rem:PDP}, in which 
% $\tilde{\lftd{B}}_{k|t}=\tilde{\lftd{B}}_{0|t+kd}$, $\tilde{\lftd{R}}_{k|t}=\tilde{\lftd{R}}_{0|t+kd}$,  $\tilde{\lftd{Q}}_{k|t}=\tilde{\lftd{Q}}_{0|t+kd}$, and $\tilde{\lftd{A}}_{k|t}=\tilde{\lftd{A}}_{0|t+kd}$ 
as defined in~\eqref{eq:btilch5},~\eqref{eq:rtilch5},~\eqref{eq:qtilch5}, and~\eqref{eq:atilch5}.
\end{hypothesis}

\begin{lemma}\label{lemma:ric_upper_bound}
Under \mc{Part 2) of Hypothesis}~\ref{asm:gamprf} on $\gamma\in\mathbb{R}_{>0}$, the lifted Riccati operator defined in~\eqref{eq:lifted_ric_0} satisfies
\begin{align*}
\tilde{\lftd{\mathcal{R}}}_{\gamma,k|t}(P) \prec \tilde{\lftd{Q}}_{k|t} + \tilde{\lftd{A}}_{\smash{k|t}}^{\trs} (\tilde{\lftd{B}}_{k|t}\tilde{\lftd{R}}_{\smash{k|t}}^{-1}\tilde{\lftd{B}}_{\smash{k|t}}^{\trs})^{-1} \tilde{\lftd{A}}_{k|t}
\end{align*}
for $P \in \mathbb{S}_{\smash{++}}^n$ \mc{and $t,k\in\mathbb{N}_0$}. 
\end{lemma}

\begin{proof}
Given $P^{-1}\in\mathbb{S}_{\smash{++}}^n$ and $\tilde{\lftd{B}}_{k|t}\tilde{\lftd{R}}_{\smash{k|t}}^{-1}\tilde{\lftd{B}}_{\smash{k|t}}^{\trs}\in\mathbb{S}_{\smash{++}}^n$, note that
$(P^{-1}+\tilde{\lftd{B}}_{k|t}\tilde{\lftd{R}}_{\smash{k|t}}^{-1}\tilde{\lftd{B}}_{\smash{k|t}}^{\trs})^{-1} \prec (\tilde{\lftd{B}}_{k|t}\tilde{\lftd{R}}_{\smash{k|t}}^{-1}\tilde{\lftd{B}}_{\smash{k|t}}^{\trs})^{-1}$.
As such, since $\tilde{\lftd{A}}_{k|t}$ is non singular by Lemma~\ref{lemma:lifted_A_nonsingular}, 
\begin{align*}
&\tilde{\lftd{A}}_{\smash{k|t}}^{\trs}(P^{-1}+\tilde{\lftd{B}}_{k|t}\tilde{\lftd{R}}_{\smash{k|t}}^{-1}\tilde{\lftd{B}}_{\smash{k|t}}^{\trs})^{-1}\tilde{\lftd{A}}_{k|t} 
\\
&\qquad \qquad \qquad
\prec \tilde{\lftd{A}}_{\smash{k|t}}^{\trs}(\tilde{\lftd{B}}_{k|t}\tilde{\lftd{R}}_{\smash{k|t}}^{-1}\tilde{\lftd{B}}_{\smash{k|t}}^{\trs})^{-1}\tilde{\lftd{A}}_{k|t}.
\end{align*}
Thus, the result holds in view of~\eqref{eq:lifted_ric_2}.
%$\tilde{\lftd{\mathcal{R}}}_{\gamma,k|t}(P) \prec \tilde{\lftd{Q}}_{k|t} + \tilde{\lftd{A}}_{\smash{k|t}}^{\trs} (\tilde{\lftd{B}}_{k|t}\tilde{\lftd{R}}_{\smash{k|t}}^{-1}\tilde{\lftd{B}}_{\smash{k|t}}^{\trs})^{-1} \tilde{\lftd{A}}_{k|t}$.
\end{proof}

\begin{lemma}\label{lemma:ric_monotone_ch5}
Under \mc{Part 2) of Hypothesis}~\ref{asm:gamprf} on $\gamma\in\mathbb{R}_{>0}$, 
\begin{align*}
P \prec X  \quad \implies \quad
\tilde{\lftd{\mathcal{R}}}_{\gamma,k|t}(P) \prec \tilde{\lftd{\mathcal{R}}}_{\gamma,k|t}(X) 
\end{align*}
for $X, P \in \mathbb{S}_{\smash{++}}^n$ \mc{and $t,k\in\mathbb{N}_0$}.
\end{lemma}
 
\begin{proof}
Given $0\prec P \prec X$ and $\tilde{\lftd{B}}_{k|t}\tilde{\lftd{R}}_{\smash{k|t}}^{-1}\tilde{\lftd{B}}_{\smash{k|t}}^{\trs}\in\mathbb{S}_{\smash{++}}^n$, note that
$0 \prec X^{-1} + \tilde{\lftd{B}}_{k|t}\tilde{\lftd{R}}_{\smash{k|t}}^{-1}\tilde{\lftd{B}}_{\smash{k|t}}^{\trs} \prec P^{-1} + \tilde{\lftd{B}}_{k|t}\tilde{\lftd{R}}_{\smash{k|t}}^{-1}\tilde{\lftd{B}}_{\smash{k|t}}^{\trs}$. Therefore, $(P^{-1} + \tilde{\lftd{B}}_{k|t}\tilde{\lftd{R}}_{\smash{k|t}}^{-1}\tilde{\lftd{B}}_{\smash{k|t}}^{\trs})^{-1}  \prec (X^{-1} + \tilde{\lftd{B}}_{k|t}\tilde{\lftd{R}}_{\smash{k|t}}^{-1}\tilde{\lftd{B}}_{\smash{k|t}}^{\trs})^{-1}$, whereby 
\begin{align*}
&\tilde{\lftd{\mathcal{R}}}_{\gamma,k|t}(P) \\
&=\tilde{\lftd{Q}}_{k|t} + \tilde{\lftd{A}}_{\smash{k|t}}^{\trs} (P^{-1} + \tilde{\lftd{B}}_{k|t}\tilde{\lftd{R}}_{\smash{k|t}}^{-1}\tilde{\lftd{B}}_{\smash{k|t}}^{\trs})^{-1} \tilde{\lftd{A}}_{k|t} \\
& \prec \tilde{\lftd{Q}}_{k|t} + \tilde{\lftd{A}}_{\smash{k|t}}^{\trs} (X^{-1} + \tilde{\lftd{B}}_{k|t}\tilde{\lftd{R}}_{\smash{k|t}}^{-1}\tilde{\lftd{B}}_{\smash{k|t}}^{\trs})^{-1} \tilde{\lftd{A}}_{k|t} = \tilde{\lftd{\mathcal{R}}}_{\gamma,k|t}(X)
\end{align*}
since $\tilde{\lftd{A}}_{k|t}$ is non-singular by
Lemma~\ref{lemma:lifted_A_nonsingular}.
\end{proof}

\begin{lemma}\label{lemma:X_bigger_P_ch5}
Under Hypothesis~\ref{asm:gamprf} on $\gamma\in\mathbb{R}_{>0}$, 
\begin{align*}
(\forall t\in\mathbb{N}_0) \quad  P_{t+1} \prec X_{t+1}
\end{align*}
with
$(X_{t+1})_{t\in\mathbb{N}_0}$ as defined in~\eqref{eq:h_inf_approx_X} given $T\in\mathbb{N}$.
\end{lemma}

\begin{proof}
By Remarks~\ref{rem:lift_equiv} and~\ref{rem:prePDP}, 
$P_{t+1+dT} = \mathcal{R}_{\gamma,t+1+dT} \circ \cdots \circ \mathcal{R}_{\gamma,t+d(T+1)}(P_{t+1+d(T+1)}) = \tilde{\lftd{\mathcal{R}}}_{\gamma,T|t+1}(P_{t+1+d(T+1)})$.
\mc{Thus}, 
\begin{align}\label{eq:compare_X_P_ch5}
P_{t+1+dT} 
&\prec  \tilde{X}_{t+1}
\end{align}
by Lemma~\ref{lemma:ric_upper_bound},
where $\tilde{X}_{t+1}$ is defined in~\eqref{eq:h_inf_terminal}.
\mc{As such},
$P_{t+1}
= \mathcal{R}_{\gamma,t+1} \circ \cdots \circ \mathcal{R}_{\gamma,t+dT}(P_{t+1+dT})
=\tilde{\lftd{\mathcal{R}}}_{\gamma,0|t+1} \circ \cdots \circ \tilde{\lftd{\mathcal{R}}}_{\gamma,T-1|t+1}(P_{t+1+dT}) 
\prec \tilde{\lftd{\mathcal{R}}}_{\gamma,0|t+1} \circ \cdots \circ \tilde{\lftd{\mathcal{R}}}_{\gamma,T-1|t+1}(\tilde{X}_{t+1})=X_{t+1}$, 
where the last equality is~\eqref{eq:h_inf_approx_X}, and
the matrix inequality follows from~\eqref{eq:compare_X_P_ch5} by Lemma~\ref{lemma:ric_monotone_ch5}.
\end{proof}

\mc{A uniform bound on} the Riemannian distance $\delta(\tilde{X}_{t+1},P_{t+1+dT})$, \mc{as per Lemma~\ref{lemma:h_inf_delta_ub}}, \mc{leads to} the same for $\delta(X_{t+1},P_{t+1})$ \mc{as per} the subsequent Lemma~\ref{lemma:h_inf_contraction_ub}. 
\begin{lemma}\label{lemma:h_inf_delta_ub}
Under \mc{Hypothesis~\ref{asm:gamprf}} on $\gamma\in\mathbb{R}_{>0}$, 
\[(\forall t\in\mathbb{N}_0)\quad \delta(\tilde{X}_{t+1},P_{t+1+dT}) \leq \overline{\delta}\] 
with
$(\tilde{X}_{t+1})_{t\in\mathbb{N}_0}\in\mathbb{S}_{\smash{+}}^n$ as defined in~\eqref{eq:h_inf_terminal} given $T\in\mathbb{N}$, where $\delta(\cdot,\cdot)$ is 
defined in~\eqref{eq:riemannian}, 
%denotes the Riemannian distance in~\eqref{eq:riemannian},
and $\overline{\delta}\in\mathbb{R}_{>0}$ is given in~\eqref{eq:deltaup}. 
% \begin{align}
% \overline{\delta} &:=
% \sqrt{n} \log \!\left(
% \sup_{t\in\mathbb{N}_0} \!
% \frac{\lambda_{\max}(\tilde{\lftd{Q}}_{0|t} + \tilde{\lftd{A}}_{\smash{0|t}}^{\trs} (\tilde{\lftd{B}}_{0|t} \tilde{\lftd{R}}_{\smash{0|t}}^{-1} \tilde{\lftd{B}}_{\smash{0|t}}^{\trs})^{-1}\tilde{\lftd{A}}_{0|t})}{\lambda_{\min}(\tilde{\lftd{Q}}_{0|t})} \right)  \label{eq:h_inf_delta_ub}
% \end{align}
\end{lemma}

\begin{proof}
From~\eqref{eq:compare_X_P_ch5},  $0\prec P_{t+1+dk} \prec \tilde{X}_{t+1}$. Therefore, using the same argument underlying~\eqref{eq:lambdamin} and~\eqref{eq:deltalambdamax},
%in the proof of Theorem~\ref{theorem:delta_requirement}, 
it follows that
$\{\lambda_1,\ldots,\lambda_n\}=\mathrm{spec}(\tilde{X}_{t+1}P_{t+1+dT}^{-1})
\subset\mathbb{R}_{>1}$, and 
\begin{align*}
\delta(\tilde{X}_{t+1},P_{t+1+dT})
&= {\textstyle \sqrt{\sum_{i=1}^n \big(\log(\lambda_i)\big)^2}}
\\
&\leq {\textstyle \sqrt{n \, \max_{i\in\{1,\ldots,n\}} \big(\log(\lambda_i)\big)^2}}
\\
&\leq \sqrt{n} \log(\lambda_{\max}(P_{t+1+dT}^{-\smash{1\!/2}}
\tilde{X}_{t+1} P_{t+1+dT}^{-\smash{1\!/2}})) \\
&\leq \sqrt{n} \log(\lambda_{\max}(\tilde{X}_{t+1}) / \lambda_{\min}(P_{t+1+dT})).
\end{align*}
From~\eqref{eq:lifted_ric_2} and Remark~\ref{rem:lift_equiv}, $\lambda_{\min}(P_{t+1+dT}) \geq \lambda_{\min}(\tilde{\lftd{Q}}_{0|t+1+dT}) \geq \inf_{t\in\mathbb{N}_0} \lambda_{\min}(\tilde{\lftd{Q}}_{0|t})$, and since
\begin{align*}
&\tilde{X}_{t+1} = \tilde{\lftd{Q}}_{0|t+1+dT} \\
&+\!\tilde{\lftd{A}}_{0|t+1+dT}^{\trs}(\tilde{\lftd{B}}_{0|t+1+dT} \tilde{\lftd{R}}_{0|t+1+dT}^{-1} \tilde{\lftd{B}}_{0|t+1+dT}^{\trs})^{-1}\!\tilde{\lftd{A}}_{0|t+1+dT},
\end{align*}
it follows that
\begin{align*}
&\delta(\tilde{X}_{t+1},P_{t+1+dT}) \\
&\qquad \!\leq\!\!
\sqrt{n} \log \!\Big(
\sup_{t\in\mathbb{N}_0} \!
\frac{\lambda_{\max}(\tilde{\lftd{Q}}_{0|t} + \tilde{\lftd{A}}_{\smash{0|t}}^{\trs} (\tilde{\lftd{B}}_{0|t} \tilde{\lftd{R}}_{\smash{0|t}}^{-1}\! \tilde{\lftd{B}}_{\smash{0|t}}^{\trs})^{-1}\!\tilde{\lftd{A}}_{0|t})}{\lambda_{\min}(\tilde{\lftd{Q}}_{0|t})} \Big).
\end{align*}
This upper bound is finite under part 2) of Hypothesis~\ref{asm:gamprf}.
\end{proof}

\begin{lemma}\label{lemma:h_inf_contraction_ub}
Under \mc{Hypothesis~\ref{asm:gamprf}} on $\gamma\in\mathbb{R}_{>0}$,
\begin{align}\label{eq:h_inf_contraction_ub}
(\forall t\in\mathbb{N}_0)~ \delta(X_{t+1},P_{t+1}) \leq (\, \overline{\rho}\,)^T \cdot \overline{\delta}
\end{align}
with
$(X_{t+1})_{t\in\mathbb{N}_0}\in\mathbb{S}_{\smash{+}}^n$ as defined in~\eqref{eq:h_inf_approx_X} given $T\in\mathbb{N}$,
and $\overline{\delta},\overline{\rho} 
\in\mathbb{R}_{>0}$ as per~\eqref{eq:deltaup} and~\eqref{eq:sup_rho_ch5},
respectively.
\end{lemma}

\begin{proof}
For $t\in\mathbb{N}_0$, $\delta(\tilde{X}_{t+1},P_{t+1+dT}) \leq \overline{\delta}$ by Lemma~\ref{lemma:h_inf_delta_ub}. 
From~\eqref{eq:h_inf_approx_X},
$X_{t+1} = \tilde{\lftd{\mathcal{R}}}_{\gamma,0|t+1} \circ \cdots \circ \tilde{\lftd{\mathcal{R}}}_{\gamma,T-1|t+1}(\tilde{X}_{t+1})$,
and by Lemma~\ref{lemma:ric_gamma_comp},
$P_{t+1} = \tilde{\lftd{\mathcal{R}}}_{\gamma,0|t+1} \circ \cdots \circ \tilde{\lftd{\mathcal{R}}}_{\gamma,T-1|t+1}(P_{t+1+dT})$.
% Since $\tilde{\lftd{R}}_{k|t}=\tilde{\lftd{R}}_{0|t+dk}$ 
% %is non-singular for all $k\in\mathbb{N}_0$, 
% and
% $\tilde{\lftd{A}}_{k|t}=\tilde{\lftd{A}}_{0|t+dk}$, 
%is also non-singular by Lemma~\ref{lemma:lifted_A_nonsingular}. 
Thus, repeated application of
Theorem~\ref{theorem:ric_gamma_contraction} yields
\begin{align*}
\delta(X_{t+1}, P_{t+1}) 
\leq (\,\overline{\rho}\,)^T \cdot \delta(\tilde{X}_{t+1},P_{t+1+dT}) 
\leq (\,\overline{\rho}\,)^T \cdot \overline{\delta}
\end{align*}
with $\overline{\rho} = \sup_{s\in\mathbb{N}_0} 1/(1+\tilde{\omega}_{0|s}) \geq 
%\sup_{k\in\mathbb{N}_0} 1/(1+\tilde{\omega}_{0|t+kd}) 
\tilde{\rho}_t= \sup_{k\in\mathbb{N}_0} 1/(1+\tilde{\omega}_{k|t})$,
$\tilde{\omega}_{k|t}:=\tilde{\epsilon}_{k|t}/\tilde{\zeta}_{k|t}\mc{=\tilde{\omega}_{0|t+kd}}$, and $\tilde{\epsilon}_{k|t},\tilde{\zeta}_{k|t}\in\mathbb{R}_{>0}$ as per
\eqref{eq:h_inf_epsilon} and \eqref{eq:h_inf_zeta}, respectively.

It remains to establish $\overline{\rho} < 1$.
%, which holds since $\inf_{s\in\mathbb{N}_0} \tilde{\omega}_{0|s} > 0$, as shown below. 
From~\eqref{eq:h_inf_zeta} and~\eqref{eq:h_inf_epsilon}, 
$\inf_{s\in\mathbb{N}_0}
\tilde{\omega}_{0|s}
\geq \epsilon_0 \big/ \zeta_0$,
where
\begin{align*}
\zeta_0&:=
\sup_{s\in\mathbb{N}_0} \zeta_{0|s}
\leq 1\big/ \inf_{s\in\mathbb{N}_0} \lambda_{\min}(\tilde{\lftd{Q}}_{0|s})
\end{align*}
because $\lambda_{\min}(\tilde{\lftd{Q}}_{0|s} + \tilde{\lftd{Q}}_{0|s}\tilde{\lftd{A}}_{\smash{0|s}}^{-1}\tilde{\lftd{B}}_{0|s}\tilde{\lftd{R}}_{\smash{0|s}}^{-1}\tilde{\lftd{B}}_{0|s}^{\trs}(\tilde{\lftd{A}}_{\smash{0|s}}^{\trs})^{-1}\tilde{\lftd{Q}}_{0|s}) \geq \lambda_{\min}(\tilde{\lftd{Q}}_{0|s})$, 
and 
\begin{align*}
\epsilon_0 &:= \inf_{s\in\mathbb{N}_0} \epsilon_{0|s}\\
&=1\big/ \sup_{s\in\mathbb{N}_0} \lambda_{\max}(\tilde{\lftd{Q}}_{0|s} + 
\tilde{\lftd{A}}_{0|s}(\tilde{\lftd{B}}_{0|s}\tilde{\lftd{R}}_{\smash{0|s}}^{-1}\tilde{\lftd{B}}_{\smash{0|s}}^{\trs})^{-1} \tilde{\lftd{A}}_{0|s}). 
\end{align*}
With the uniform bounds in part 2) of Hypothesis~\ref{asm:gamprf}, $\zeta_0<+\infty$ and $\epsilon_0>0$, whereby $\epsilon_0/\zeta_0 >0$.
\end{proof}

%\begin{proof}
{\noindent\hspace{1em}{\bf\itshape Proof of 
Theorem~\ref{theorem:h_inf_T_lb}:}
 \mc{Given} the strict inequality~\eqref{eq:preview_lb}, \mc{with $\mc{\eta}$ as defined in~\eqref{eq:alpha_ch5},} direct calculation yields
\begin{align*}
(\exists \varepsilon\in\mathbb{R}_{>0}) \quad 
(\overline{\rho})^T \cdot \overline{\delta} \leq 
\log((\mc{\eta}-\varepsilon)\cdot\underline{\kappa}  + 1).
\end{align*}
 By Lemma~\ref{lemma:h_inf_contraction_ub}, $(\forall t\in\mathbb{N}_0) \, \delta(X_{t+1}, P_{t+1}) \leq (\, \overline{\rho}\,)^T \cdot \overline{\delta}$. Thus, $(\exists \varepsilon\in\mathbb{R}_{>0})~(\forall t \in \mathbb{N}_0)$
\begin{align*}
\delta(X_{t+1}, P_{t+1}) &\leq \log((\mc{\eta}-\varepsilon)\cdot \underline{\kappa} + 1)\\
&\leq \log((\mc{\eta}-\varepsilon)\cdot \underline{\lambda} + 1),
\end{align*}
since $\inf_{t\in\mathbb{N}_0} \lambda_{\min}(\tilde{\lftd{Q}}_{0|t}) =:\underline{\kappa}\leq \underline{\lambda} := \inf_{t\in\mathbb{N}_0} \lambda_{\min}(P_t)$ in view of Remark~\ref{rem:lift_equiv} and~\eqref{eq:lifted_ric_2}.
%; see Remark~\ref{rem:h_inf_Pk_PD}.
Further, $P_{t+1} \prec X_{t+1}$ by Lemma~\ref{lemma:X_bigger_P_ch5}. 
Therefore, the sufficient condition in~\eqref{eq:X_requirement_final} is satisfied, and Theorem~\ref{theorem:delta_requirement} applies to yield Theorem~\ref{theorem:h_inf_T_lb}.
\hspace*{\fill}\QEDopen

\section{Numerical Example} \label{sec:numex}

\mc{A periodic numerical example is presented here to explore aspects of Theorem~\ref{theorem:h_inf_T_lb}. As far as the authors are aware, this main result admits the only practical synthesis of a state feedback control policy guaranteed to achieve a worst-case $\ell_2$ gain specification in a general time-varying setting. In the special case of a periodic system, the relevant constants $\underline{\kappa}$, $\overline{\delta}$, and $\overline{\rho}$ can be determined exactly, which enables assessment of the analytical conservativeness noted in Remark~\ref{rem:bounds}, and elaborated below in the context of the linearization of a nonlinear system along a periodic trajectory.} 

With sampling interval $h\in\mathbb{R}_{>0}$, Euler discretization of the unicycle kinematics yields the state space model
\begin{align*}
\begin{bmatrix}
x_{t+1} \\ y_{t+1} \\ \mc{\psi}_{t+1}
\end{bmatrix}
= \begin{bmatrix}
x_t + v_t \cos(\mc{\psi}_t) h \\
y_t + v_t \sin(\mc{\psi}_t) h \\
\mc{\psi}_t + \mc{r}_t h
\end{bmatrix},
\end{align*}
where $(x_t, y_t)$ is the position of the robot in the plane, $\mc{\psi}_t$ is the yaw angle, $v_t$ is the velocity along the yaw angle, and $\mc{r_t}$ is the yaw rate, at continuous times $t\cdot h\in\mathbb{R}_{\geq 0}$, for $t\in\mathbb{N}_0$. 

Consider the nominal lemniscate state trajectory given by
\begin{align*}
&x_t^* = \frac{a \cos(k)}{1+(\sin(k))^2}, \quad y_t^* = \frac{a \sin(k)\cos(k)}{1+(\sin(k))^2},\\
&\qquad\quad \mc{\psi}_t^* = \arctan\left( \frac{y_{t+1}^*-y_t^*}{x_{t+1}^*-x_t^*} \right),
\end{align*}
for $t\in\mathbb{N}_0$,
where
$k = \frac{2\pi}{N} (t \mod ~N)$,
%\mathrm{rem}(t,N)$,
$a=1$ is the half-width of the lemniscate, and $N \in\mathbb{N}$ is the period. The nominal trajectory in the $xy$-plane is shown in Figure~\ref{fig:nominal_trajectory} for $N=400$. The corresponding nominal inputs are given by
\begin{align*}
v_t^* &= \sqrt{(x_{t+1}^*-x_{t}^*)^2 + (y_{t+1}^*-y_{t}^*)^2} \Big/ h, \\
\mc{r}_t^* &= (\mc{\psi}_{t+1}^*-\mc{\psi}_t^*)/h.
\end{align*}
Linearization along the nominal trajectory yields the time-varying discrete time model
\begin{align*}
\begin{bmatrix}
x_{t+1} - x_{t+1}^* \\
y_{t+1} - y_{t+1}^* \\
\mc{\psi}_{t+1} - \mc{\psi}_{t+1}^*
\end{bmatrix}
&=
\begin{bmatrix}
1 & 0 & -v_t^* \sin(\mc{\psi}_t^*) h \\
0 & 1 & v_t^* \cos(\mc{\psi}_t^*) h \\
0 & 0 & 1
\end{bmatrix}
\begin{bmatrix}
x_t-x_t^* \\ y_t-y_t^* \\ \mc{\psi}_t-\mc{\psi}_t^*
\end{bmatrix}
\\
&\quad + 
\begin{bmatrix}
\cos(\mc{\psi}_t^*)h & 0 \\
\sin(\mc{\psi}_t^*)h & 0 \\
0 & h
\end{bmatrix}
\begin{bmatrix}
v_t-v_t^* \\ \mc{r}_t-\mc{r}_t^*
\end{bmatrix} + h \, w_t,
\end{align*}
in the form~\eqref{eq:ltv_sys_w}
with disturbance input $h\, w_t$. 

\begin{figure}[t!]
    \centering
    \includegraphics[width=0.98\columnwidth]{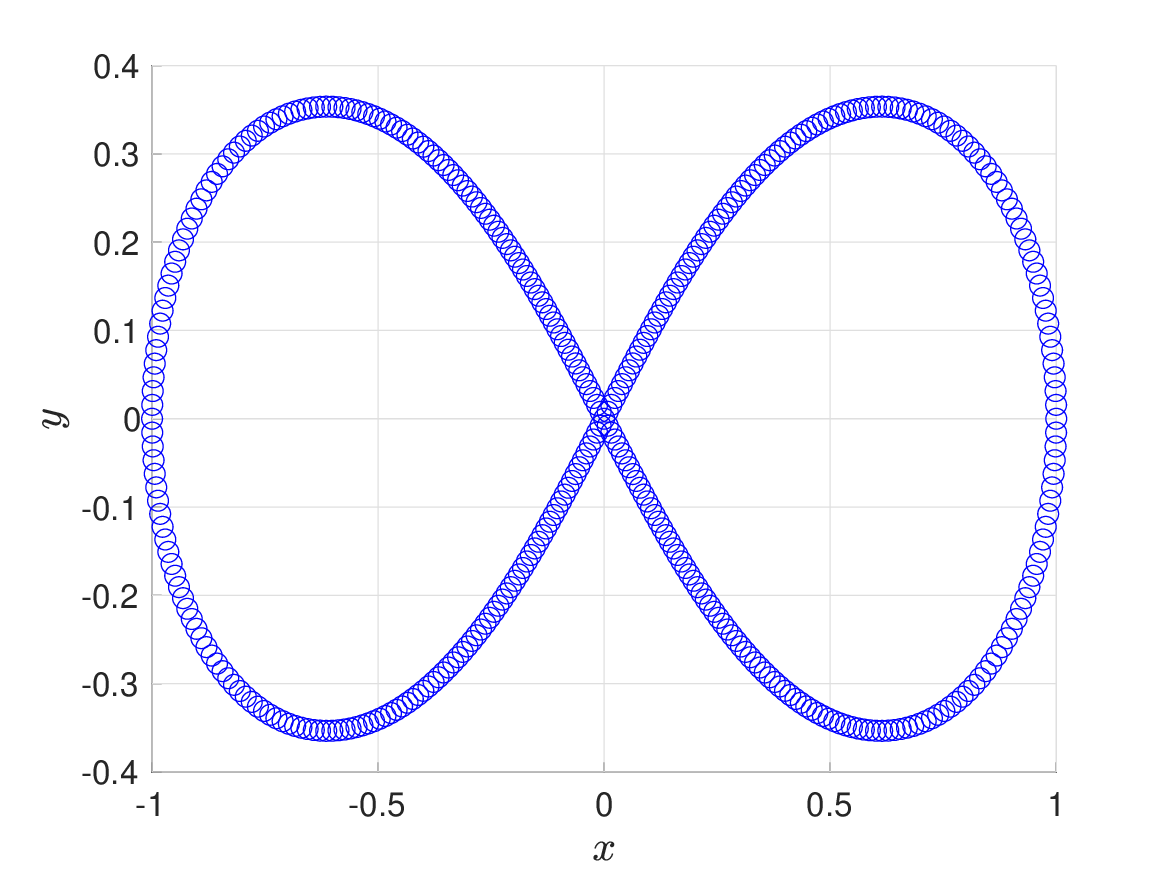}
    \caption{Nominal trajectory.}
    \label{fig:nominal_trajectory}
\end{figure}

For the performance output $z_t$ in~\eqref{eq:output_z_def} with
\begin{align*}
    Q_{t} = \begin{bmatrix} 2 &  0 &  0\\
            0 & 2 &  0 \\
            0 & 0 & 0.2 \end{bmatrix}
            \quad \text{ and }\quad
    R_{t} = \begin{bmatrix} 0.1 & 0 \\
            0 & 0.01 \end{bmatrix},
\end{align*}
it can be verified
with sampling interval $h=0.05$, and baseline $\ell_2$ gain bound $\gamma=125$, that parts 1) and 2) of the hypothesis in Theorem~\ref{theorem:h_inf_T_lb} hold for all lifting steps $d\in[10:40]$, set in accordance with Assumption~\ref{asm:uni_obs_ctr}. This amounts to checking the finite number of conditions across one period of the problem data for part 2), and exploiting the periodic structure of the Riccati recursion~\eqref{eq:RiccatiRecursion} as shown in, e.g.,~\cite{hench1994numerical,bittanti2009periodic}. In particular, it can be established that there exists periodic $(P_t)_{t\in\mathbb{N}_0}\subset\mathbb{S}_{\smash{++}}^3$ such that~\eqref{eq:sc_P_ub} and~\eqref{eq:RiccatiRecursion} both hold, and the infinite-preview control policy~\eqref{eq:h_inf_causal_policy} achieves $(\forall w\in\ell_2)~\|z\|_2 \leq \gamma \|h\, w\|_2 = 6.25\|w\|_2$. 

\begin{figure}[t!]
    \centering
    \includegraphics[width=0.98\columnwidth]{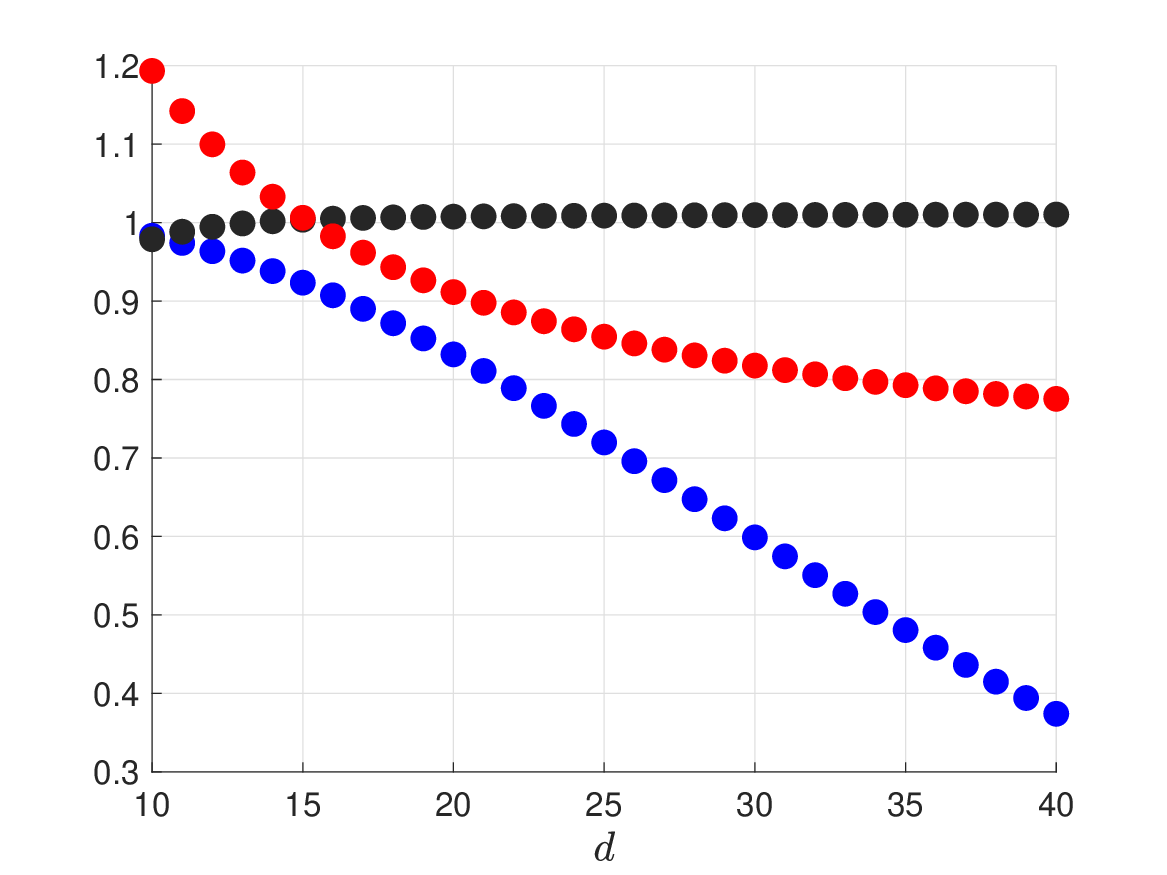}
    \caption{$\underline{\kappa}$ (black), $0.1\overline{\delta}$ (red), and $\overline{\rho}$ (blue), in Theorem~\ref{theorem:h_inf_T_lb} for lifting steps $d\in[10:40]$.}
    \label{fig:Thm2params}
\end{figure}

\begin{figure}[!t]
\centering
\includegraphics[width=0.98\columnwidth]{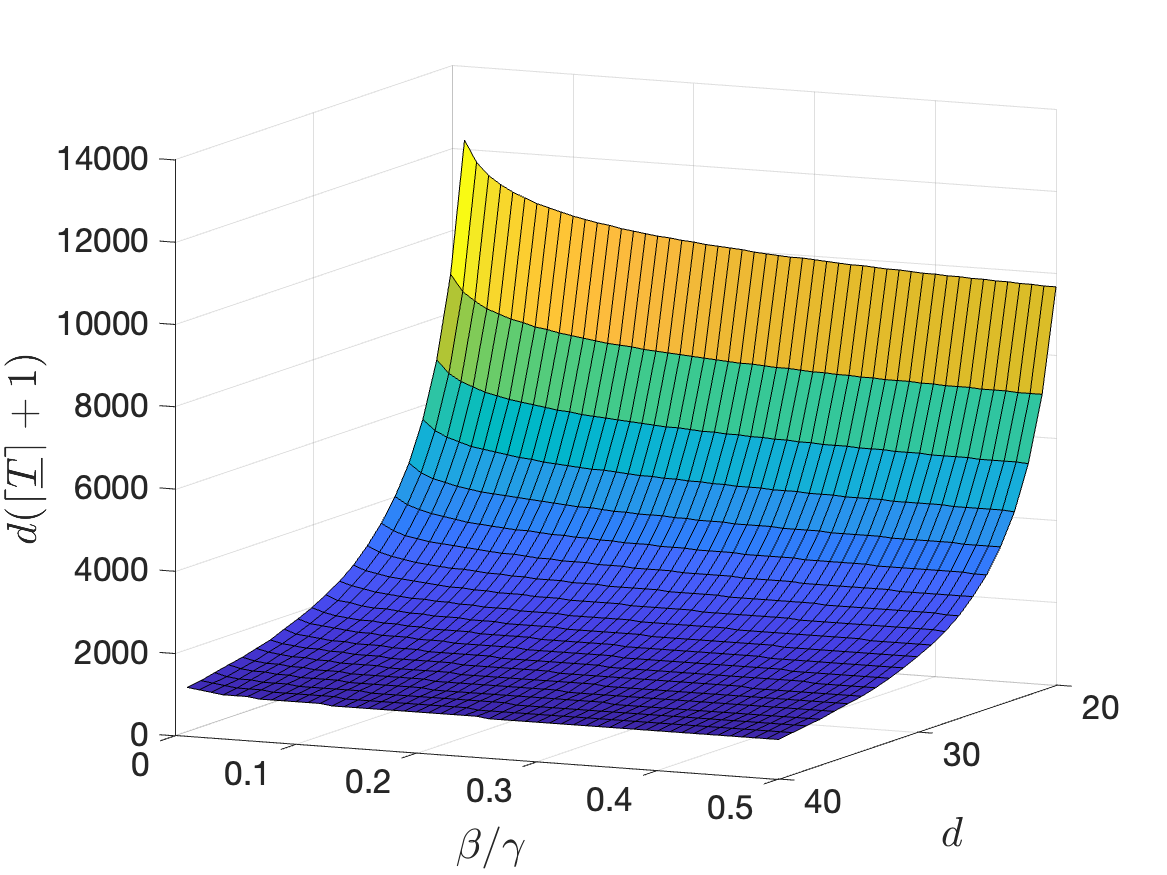}
\caption{Preview steps $d(\lceil\underline{T}\rceil+1)$ sufficient \mc{according to Theorem~\ref{theorem:h_inf_T_lb}} for \mc{different} performance loss bounds $\beta$ \mc{relative to baseline gain bound $\gamma=125$ ($1-50\%$)}, and lifting steps $d\in[20:40]$.}
\label{fig:dTvBetaD}
\end{figure}

\begin{figure}[t!]
    \centering
    \includegraphics[width=0.98\columnwidth]{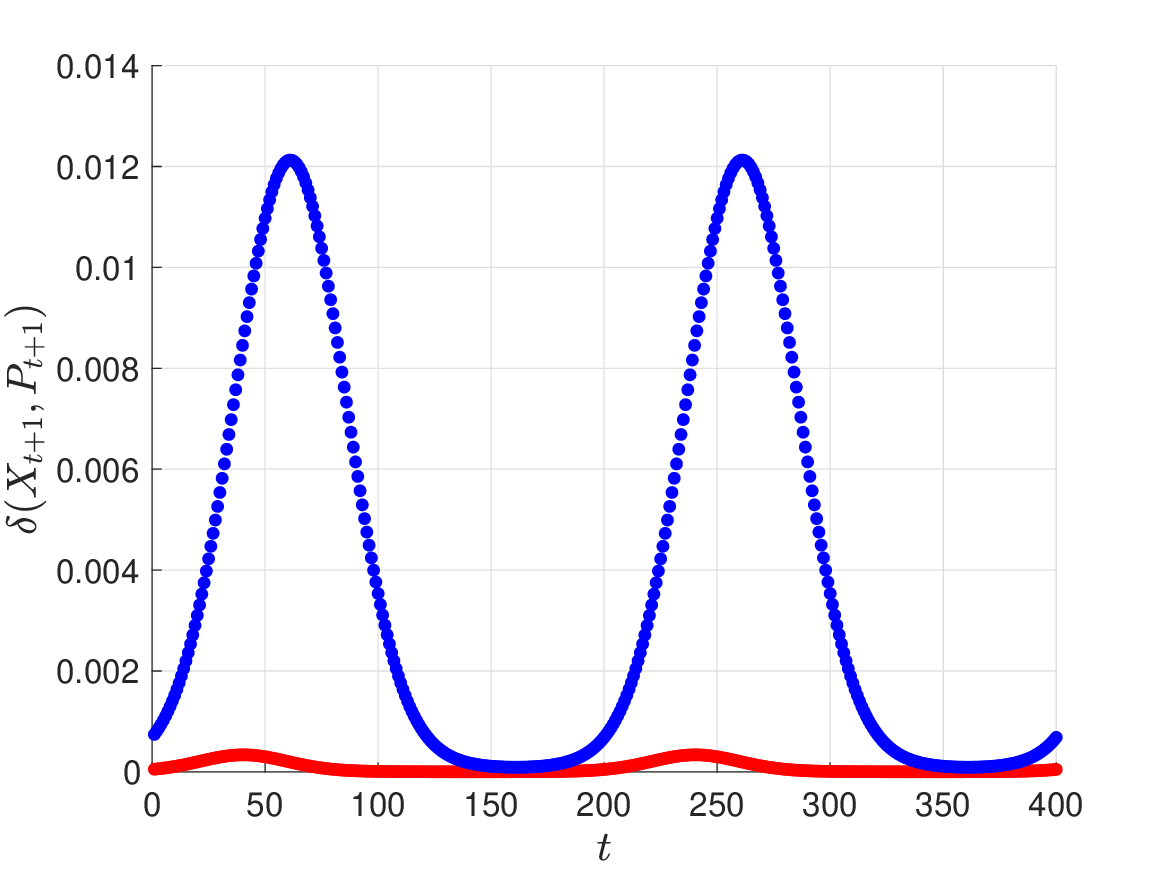}
    \caption{$\delta(X_{t+1},P_{t+1})$ for: (i) $T=1$, $d=40$ $\mc{\implies} \overline{\rho}^T\overline{\delta}=2.90$ (blue); and (ii) $T=2$, $d=40$ $\mc{\implies} \overline{\rho}^T\overline{\delta}=1.08$ (red).}
    \label{fig:deltaXP}
\end{figure}

In Figure~\ref{fig:Thm2params}, it is shown how $\underline{\kappa}$, $\overline{\delta}$, and $\overline{\rho}$ in Theorem~\ref{theorem:h_inf_T_lb} depend on the number of lifting steps $d\in[10:40]$; \mc{see Remark~\ref{rem:ddep}}. In the periodic setting of this example it is possible to evaluate these quantities exactly; one period suffices to determine the $\sup$ and $\inf$ in~\eqref{eq:kappalo}--\eqref{eq:omtil}. Increasing $d$ decreases both the approximate Riccati solution error bound $\overline{\delta}$ from Lemma~\ref{lemma:h_inf_delta_ub}, and the  Riccati contraction rate bound $\overline{\rho}$ in the proof of Lemma~\ref{lemma:h_inf_contraction_ub}, whereas $\underline{\kappa}$ remains roughly constant. 

The \mc{number} $d(\lceil\underline{T}\rceil+1)$ of preview steps that is sufficient for the control 
policy~\eqref{eq:h_inf_approx_policy}, with $(X_{t+1})_{t\in\mathbb{N}_0}$ as per~\eqref{eq:h_inf_approx_X}, to achieve the performance loss bound $\beta$, relative to \mc{the} baseline gain bound $\gamma=125$, is shown in Figure~\ref{fig:dTvBetaD} for $d\in[20:40]$. When $d=40$, this estimate of \mc{a sufficient number of} preview steps lies between only $2$ and $3$ system periods ($800-1200$ steps) for all $\beta\in[0.01\gamma,0.5\gamma]$, i.e., $1-50\%$ performance loss. The measured approximation error $\delta(X_{t+1},P_{t+1})\leq(\overline{\rho})^T\overline{\delta}$ (see Lemma~\ref{lemma:h_inf_contraction_ub}) is shown in Figure~\ref{fig:deltaXP} for $T=1$ and $T=2$, which corresponds to \mc{just} $80$ and $120$ preview horizon steps, respectively. In both cases, the measured Riemannian distance is two to three orders of magnitude smaller than the bound used to prove Theorem~\ref{theorem:h_inf_T_lb}.

\section{Conclusion} \label{sec:conc}

A finite \mc{receding-horizon} preview approximation of \mc{the standard} infinite-preview state feedback $\ell_2$ gain control policy is considered \mc{for linear time-varying systems}. Under uniform controllability and observability, and given a suitably large baseline gain bound, the strict contraction of lifted Riccati operators is exploited to construct a feedback gain at each step that depends on a finite-horizon preview of the problem data. The approximation achieves a \mc{prescribed} infinite-horizon performance loss bound. \mc{As such,} the proposed approach constitutes a practical controller synthesis \mc{method} with a \mc{prescribed} closed-loop $\ell_2$ gain performance \mc{guarantee}. 

As future work, it would be of interest to extend the \mc{proposed} synthesis to dynamic output feedback controllers for time-varying systems. \mc{In this case, the standard} policy structure is more \mc{complicated}, and the main challenge lies in \mc{effectively} accounting for the coupling between the two Riccati recursions that arise in the analysis of performance loss. \mc{Relaxing controllability and observability to stabilizability and detectability could also be considered.} Another direction \mc{is} to \mc{apply} the approach in formulating terminal ingredients for receding horizon schemes with hard input and state constraints, along the lines of~\cite{goulart2009control,karapetyan2022regret,martin2025guarantees} for time-invariant systems. \mc{Finally}, \mc{investigating non-stationary} structures \mc{that constrain variation in the model} data to enable tighter estimates of analysis bounds would also be of interest.

\appendix

\begin{lemma}\label{lemma:lifted_A_nonsingular}
Given $t,k\in\mathbb{N}_0$, if $\tilde{\lftd{R}}_{k|t}$ in~\eqref{eq:rtilch5} is non-singular, then $\tilde{\lftd{A}}_{k|t}$ in~\eqref{eq:atilch5} is non-singular. Further,    $\tilde{\lftd{A}}_{\smash{k|t}}^{-1}$ is uniformly bounded if $\tilde{\lftd{R}}_{\smash{k|t}}^{-1}$ is uniformly bounded. 
\end{lemma}

\begin{proof}
Consider the matrix
\begin{align}\label{eq:Hhat}
\lftd{H}_{k|t}:=\begin{bmatrix}
\lftd{A}_{k|t} & \begin{bmatrix} \lftd{B}_{k|t} & \lftd{F}_{k|t} \end{bmatrix} \\
\begin{bmatrix} \lftd{D}_{\smash{k|t}}^{\trs} \\ \lftd{E}_{\smash{k|t}}^{\trs} \end{bmatrix} \lftd{C}_{k|t} & \tilde{\lftd{R}}_{k|t}
\end{bmatrix},
\end{align}
where $\lftd{A}_{k|t}$, $\lftd{B}_{k|t}$, $\lftd{F}_{k|t}$, $\lftd{C}_{k|t}$, $\lftd{D}_{k|t}$, and $\lftd{E}_{k|t}$, are as defined in Lemmas~\ref{lemma:lifted_x} and~\ref{lemma:lifted_z}.
With $\tilde{\lftd{R}}_{k|t}$ non-singular, 
the corresponding Schur complement of $\lftd{H}_{k|t}$ is equal to $\tilde{\lftd{A}}_{k|t}$. As such, $\lftd{H}_{k|t}$ is non-singular if and only if $\tilde{\lftd{A}}_{k|t}$ is non-singular. By Assumption~\ref{asm:boundeddata},
$\lftd{A}_{k|t} =  \varPhi_{t+dt+d,t+dk}$ is also invertible. Thus, the alternative Schur complement of $\lftd{H}_{k|t}$ in~\eqref{eq:Hhat} is given by
\begin{align}\label{eq:proof_a_nonsingular_schur}
\lftd{S}_{k|t}
&= \tilde{\lftd{R}}_{k|t} - \begin{bmatrix} \lftd{D}_{\smash{k|t}}^{\trs} \\ \lftd{E}_{\smash{k|t}}^{\trs} \end{bmatrix} \lftd{C}_{k|t} \lftd{A}_{\smash{k|t}}^{-1} \begin{bmatrix} \lftd{B}_{k|t} & \lftd{F}_{k|t} \end{bmatrix} 
=\mc{\varUpsilon_{k|t}} + \lftd{U}_{k|t},
\end{align}
where 
\begin{align*}
\mc{\varUpsilon_{k|t}} := \begin{bmatrix}
\lftd{R}_{k|t} & 0_{md,nd} \\ 0_{nd,md} & -\gamma^2I_{nd}
\end{bmatrix}
\end{align*}
and
\begin{align*}
&\lftd{U}_{k|t} := \begin{bmatrix} \lftd{D}_{\smash{k|t}}^{\trs} \\ \lftd{E}_{\smash{k|t}}^{\trs} \end{bmatrix} 
\begin{bmatrix} \lftd{D}_{k|t} & \lftd{E}_{k|t} \end{bmatrix} 
\!-\! \begin{bmatrix} \lftd{D}_{\smash{k|t}}^{\trs} \\ \lftd{E}_{\smash{k|t}}^{\trs} \end{bmatrix} \lftd{C}_{k|t} \lftd{A}_{\smash{k|t}}^{-1} \begin{bmatrix} \lftd{B}_{k|t} & \lftd{F}_{k|t} \end{bmatrix} \\
&=
\begin{bmatrix}
    \mc{\varXi_{k|t}}^\trs \\
    \begin{bmatrix} 0_{nd,n} & I_{nd} \end{bmatrix}
\end{bmatrix}
(\mc{\varLambda_{\smash{k|t}}^{-1}})^\trs
\mc{\varGamma_{k|t}}^\trs
\\
& \times
\mc{\varGamma_{k|t}}
(\mc{\varLambda_{\smash{k|t}}^{-1}} -
\mc{\varLambda_{\smash{k|t}}^{-1}}
\begin{bmatrix}
I_n \\ 0_{nd,n}
\end{bmatrix}
\lftd{A}_{\smash{k|t}}^{-1}
\begin{bmatrix}
0_{n,nd} & I_n
\end{bmatrix}
\mc{\varLambda_{\smash{k|t}}^{-1}}
) \\
&\times
\begin{bmatrix}
    \mc{\varXi_{k|t}}^\trs \\
    \begin{bmatrix} 0_{nd,n} & I_{nd} \end{bmatrix}
\end{bmatrix}^\trs
\end{align*}
with $\mc{\varLambda_{k|t}}$, $\mc{\varXi_{k|t}}$, and $\mc{\varGamma_{k|t}}$ as in~\eqref{eq:Ahat_ch5},~\eqref{eq:Bhat_ch5}, and~\eqref{eq:Cbreve}, respectively.
Therefore, bounded invertibility of $\lftd{S}_{k|t}=\mc{\varUpsilon_{k|t}} + \lftd{U}_{k|t}$ is equivalent to bounded invertibility of $\tilde{\lftd{A}}_{k|t}$. 

Since
\[
\mc{\varLambda_{\smash{k|t}}^{-1}}
=\begin{bmatrix}
    \lftd{G}_{k|t} & \varDelta_{k|t} \\ \lftd{A}_{k|t} & \lftd{F}_{k|t}
\end{bmatrix},
\]
where 
\begin{align*}
\lftd{G}_{k|t} = \begin{bmatrix}
I_n &
\varPhi_{t+dk+1,t+dk}^\trs & \cdots & \varPhi_{t+d(k+1)-1,t+dk}^\trs
\end{bmatrix}^\trs,   
\end{align*}
and $\varDelta_{k|t}\in\mathbb{R}^{nd\times nd}$ is block lower triangular with $n\times n$ zero blocks on the diagonal, it follows by direct calculation that 
\begin{align*}
& (\mc{\varLambda_{\smash{k|t}}^{-1}} -
\mc{\varLambda_{\smash{k|t}}^{-1}}
\begin{bmatrix}
I_n \\ 0_{nd,n}
\end{bmatrix}
\lftd{A}_{\smash{k|t}}^{-1}
\begin{bmatrix}
0_{n,nd} & I_n
\end{bmatrix}
\mc{\varLambda_{\smash{k|t}}^{-1}}
) \\
&\qquad\qquad\qquad\qquad
=\begin{bmatrix}
    0_{nd,n} & \varDelta_{k|t} - \lftd{G}_{k|t}\lftd{A}_{\smash{k|t}}^{-1}\lftd{F}_{k|t}\\
    0_{n,n} & 0_{n,nd}
\end{bmatrix},  
\end{align*}
where $\varDelta_{k|t}-\lftd{G}_{k|t}\lftd{A}_{\smash{k|t}}^{-1}\lftd{F}_{k|t} = \varDelta_{k|t} - \lftd{G}_{k|t} \varPhi_{t+d(k+1),t+dk}^{-1} \lftd{F}_{k|t}$ is block upper triangular. Similarly, 
\begin{align*}
&\begin{bmatrix}
    \mc{\varXi_{k|t}}^\trs \\
    \begin{bmatrix} 0_{nd,n} & I_{nd} \end{bmatrix}
\end{bmatrix}
(\mc{\varLambda_{\smash{k|t}}^{-1}})^\trs
\mc{\varGamma_{k|t}}^\trs
\\
&= 
\begin{bmatrix}
    \begin{bmatrix} 0_{nd,n} & \mathrm{diag}(B_{t+dk}^\trs,\ldots,B_{t+d(k+1)-1}^\trs) \end{bmatrix} \\
    \begin{bmatrix} 0_{nd,n} & I_{nd} \end{bmatrix}
\end{bmatrix}\\
&
\quad \times
\begin{bmatrix}
    \lftd{G}_{k|t}^\trs & \lftd{A}_{k|t}^\trs \\
    \varDelta_{k|t}^\trs & \lftd{F}_{k|t}^\trs
\end{bmatrix}
\begin{bmatrix}
    \mathrm{diag}((Q_{t+kd}^{\smash 1\!/2})^\trs,\ldots,(Q_{t+d(k+1)-1}^{\smash 1\!/2})^\trs) \\ 0_{nd,n}
\end{bmatrix}\\
&=
\mathrm{diag}(B_{t+dk}^\trs,\ldots,B_{t+d(k+1)-1}^\trs)\\
&\quad \times
\varDelta_{k|t}^\trs
~\mathrm{diag}((Q_{t+kd}^{\smash 1\!/2}),\ldots,(Q_{t+d(k+1)-1}^{\smash 1\!/2}))
\end{align*}
is block upper triangular with $n\times n$ zero blocks along the diagonal. As such, $\lftd{U}_{k|t}$ is block upper triangular, with all zero diagonal blocks, and thus,
\begin{align*}
\mathrm{spec}(\lftd{S}_{k|t}) &= \mathrm{spec}(\mc{\varUpsilon_{k|t}} + \lftd{U}_{k|t}) \\
&= \mathrm{spec}(\mc{\varUpsilon_{k|t}}) 
=: \{\lambda_1(\mc{\varUpsilon_{k|t}}), \ldots, \lambda_{(m+n)d}(\mc{\varUpsilon_{k|t}}) \},
\end{align*}
with \mc{real valued} $\lambda_1(\mc{\varUpsilon_{k|t}})\geq \cdots \geq \lambda_{(m+n)d}(\mc{\varUpsilon_{k|t}})$, and
\[c:=\min\{|\lambda_{i}(\mc{\varUpsilon_{k|t}})|~:~i\in\{1,\ldots,(m+n)d\}\} 
\in\mathbb{R}_{>0}.\] Application of Weyl's inequality for singular values~\cite[Thm.~3.3.2]{horn1991topics} gives
\begin{align*}
    0 < c^{(m+n) d} &\leq |\!\!\prod_{i=1}^{(m+n)d} \!\!\! \lambda_i(\mc{\varUpsilon_{k|t}})\,| 
    \\
    &=
    |\!\!\prod_{i=1}^{(m+n)d} \!\!\! \lambda_i(\lftd{S}_{k|t})\,| =
    \prod_{i=1}^{(m+n)d} \!\!\! \sqrt{\lambda_i(\lftd{S}_{k|t}^\trs \lftd{S}_{k|t})}.
\end{align*}
Hence, the smallest eigenvalue of $\lftd{S}_{k|t}^\trs \lftd{S}_{k|t}$ satisfies
\begin{align*}
    \sqrt{\lambda_{(m+n)d}(\lftd{S}_{k|t}^\trs \lftd{S}_{k|t})}
    &\geq 
    \frac{c^{(m+n) d}}{ \left( \sqrt{\lambda_{1}({S}_{k|t}^\prime {S}_{k|t}) } \right)^{(m+n) d-1} },
    % \frac{c^{(m+n) d}}{
    % \left(\sqrt{\lambda_{1}(\lftd{S}_{k|t}^\trs \lftd{S}_{k|t} )}\right)^{((m+n) d-1)}},
    %\\
    %&\geq c \in\mathbb{R}_{>0},
\end{align*}
which yields a uniform positive lower bound. In particular, from~\eqref{eq:proof_a_nonsingular_schur}, $\lambda_{1}(\lftd{S}_{k|t}^\trs \lftd{S}_{k|t})$ is uniformly bounded by Assumption~\ref{asm:boundeddata}. Therefore, $\lftd{S}_{k|t}$ is non-singular with uniformly bounded inverse, whereby the same holds for $\tilde{\lftd{A}}_{k|t}$, as claimed. 
\end{proof}

\section*{References}
\bibliography{Bibliography}
\bibliographystyle{ieeetr}

\end{document}